\documentclass[12pt]{amsart}
\usepackage[dvips]{graphics}
\usepackage[latin1]{inputenc}
\usepackage{amssymb, amsmath,   mathrsfs}
\input xy
\xyoption{all}
\DeclareMathAlphabet{\mathpzc}{OT1}{pzc}{m}{it}

\newtheorem{theorem}{Theorem}[section]

\newtheorem{proposition}[theorem]{Proposition}

\newtheorem{lemma}[theorem]{Lemma}

\theoremstyle{definition}

\theoremstyle{remark}
\newtheorem{remark}[theorem]{Remark}

\def\varle{\leqslant}
\def\varge{\geqslant}

\newcommand{\CA}{{\mathcal A}}
\newcommand{\CB}{{\mathcal B}}
\newcommand{\CC}{{\mathcal C}}
\newcommand{\CE}{{\mathcal E}}

\newcommand{\CG}{{\mathcal G}}

\newcommand{\CI}{{\mathcal I}}
\newcommand{\CJ}{{\mathcal J}}

\newcommand{\CM}{{\mathcal M}}

\newcommand{\CO}{{\mathcal O}}
\newcommand{\CP}{{\mathcal P}}

\newcommand{\CS}{{\mathcal S}}
\newcommand{\CT}{{\mathcal T}}

\newcommand{\CV}{{\mathcal V}}
\newcommand{\CW}{{\mathcal W}}

\newcommand{\Sh}{{\mathcal{SH}}}

\newcommand{\CZ}{{\mathcal Z}}

\newcommand{\SA}{{\mathscr A}}
\newcommand{\SB}{{\mathscr B}}
\newcommand{\SC}{{\mathscr C}}

\newcommand{\SL}{{\mathscr L}}
\newcommand{\SM}{{\mathscr M}}
\newcommand{\SN}{{\mathscr N}}

\newcommand{\SP}{{\mathscr P}}

\newcommand{\SV}{{\mathscr V}}

\newcommand{\fh}{{{\mathfrak h}}} 
 
\newcommand{\fp}{{{\mathfrak p}}} 
\newcommand{\fg}{{{\mathfrak g}}} 
\newcommand{\fb}{{{\mathfrak b}}}

\newcommand{\fhd}{\fh^\star}

\newcommand{\DC}{{\mathbb C}}
\newcommand{\DP}{{\mathbb P}}

\newcommand{\DZ}{{\mathbb Z}}

\newcommand{\DN}{{\mathbb N}}

\newcommand{\DV}{{\mathbb V}}

\newcommand{\DW}{{\mathbb W}}

\newcommand{\height}{{\operatorname{ht}}}

\newcommand{\im}{{\operatorname{im}}}
\newcommand{\Spec}{{\operatorname{Spec}}}
\newcommand{\End}{{\operatorname{End}}}

\newcommand{\Hom}{{\operatorname{Hom}}}
\newcommand{\id}{{\operatorname{id}}}

\newcommand{\supp}{{\operatorname{supp}}}
\newcommand{\catmod}{{\operatorname{-mod}}}

\newcommand{\inj}{{\hookrightarrow}}
\newcommand{\Loc}{{\SL}}

\newcommand{\Stab}{{\operatorname{Stab}}}

\newcommand{\linie}{{\,\text{---\!\!\!---}\,}}
\newcommand{\llinie}{{\text{---\!\!\!---\!\!\!---}}}
\newcommand{\lllinie}{{\,\text{---\!\!\!\!---\!\!\!\!---\!\!\!\!---\!\!\!\!---\!\!\!\!---\!\!\!\!---}\,}}

\begin{document}

\pagenumbering{arabic}
\title[Sheaves on moment graphs]{Sheaves on moment graphs and a localization of Verma flags}
\author[]{Peter Fiebig${}^\ast$}
\thanks{${}^\ast$ supported by the Deutsche Forschungsgemeinschaft}
\address{Mathematisches Institut, Universit{ä}t Freiburg, 79104 Freiburg, Germany}
\email{peter.fiebig@math.uni-freiburg.de}
\begin{abstract} To any moment  graph $\CG$ we assign a subcategory $\CV$ of the category of sheaves on $\CG$ together with an exact structure. We show that in the case that the graph is associated to a non-critical block of the equivariant category $\CO$ over a symmetrizable Kac--Moody algebra, $\CV$ is equivalent (as an exact category) to the subcategory of modules that admit a Verma flag. The projective modules correspond under this equivalence to the intersection cohomology sheaves on the graph.
\end{abstract}
\maketitle

\section{Introduction}
Duality constructions between geometric and representation-theoretic situations sometimes become quite transparent when interpreted in a suitable combinatorial context. Finding a combinatorial interpretation can even be crucial for establishing a duality. 

In this paper we consider the combinatorics underlying the Koszul duality between the representation category $\CO$ over a symmetrizable Kac-Moody algebra and the geometry of  Schubert varieties. In \cite{Soe1} Soergel showed, in the case of semi-simple finite-dimensional complex Lie algebras, that the projective objects in $\CO$ and the intersection cohomologies of Schubert varieties in the Langlands-dual flag variety share the same combinatorial description in terms of the algebra of coinvariants of the associated Weyl group. The proof uses translation functors and it is very difficult to describe the resulting objects (this is the main concern of combinatorial Kazhdan-Lusztig theory, cf.~ \cite{Soe3} and the references therein). 

In \cite{BMP} Braden and MacPherson gave another description of the (equivariant) intersection cohomologies on complex projective varieties with an equivariantly formal torus action in terms of a canonical sheaf on the associated moment graph.  In the case of Schubert varieties their description is essentially equivalent to Soergel's, but it has the advantage of avoiding the translation functors. In this paper we establish the Koszul-dual result by localizing the projective objects in an equivariant version of $\CO$ on the moment graph. 

Let $\CM\subset\CO$ be the subcategory of modules  admitting a Verma flag. We construct an equivalence between an equivariant version of $\CM$ and a subcategory $\CV$ of the category of sheaves on the associated moment graph. Both categories carry an exact structure and our equivalence is exact. We show that the Braden-MacPherson canonical sheaves are the indecomposable projective objects in $\CV$. Therefore they correspond to the indecomposable projective objects in $\CM$. Our approach works in the infinite dimensional context of symmetrizable Kac--Moody algebras when we exclude the blocks in the critical hyperplanes. 

We describe the ideas in this paper in some more detail. 
Let $V$ be a vector space over a field $k$. A $V$-moment graph $\CG$ is a (not necessarily finite) partially ordered graph whose edges are labelled by one-dimensional subspaces of $V$. Moment graphs occur in the following situations:
\begin{itemize}
\item As a combinatorial model for the topology of a complex equivariantly formal variety (cf.~ \cite{GKM,BMP}). The moment graph encodes the structure of the 0- and 1-dimensional orbits of a torus action.
\item As a combinatorial model for the extension structure of Verma modules in a non-critical block of the category $\CO$ over a symmetrizable Kac-Moody algebra. 
\end{itemize}

Let $S=S(V)$ be the graded symmetric algebra over $V$.
To a moment graph $\CG$ one associates the graded algebra $\CZ=\CZ(\CG)$ of tuples $(z_x)\in\prod S$, indexed by the vertices of $\CG$, such that $z_x-z_y\in S\cdot L$ if $x$ and $y$ are linked by an edge that is labelled by $L\subset V$. We call $\CZ$ the {\em structure algebra} of $\CG$. Note that it does not depend on the order on $\CG$. If $\CG$ is associated to an equivariantly formal variety, $\CZ$ is naturally isomorphic to its equivariant cohomology. If $\CG$ is associated to a block of category $\CO$, $\CZ$ is isomorphic to the equivariant center of the block.

The structure algebra is the algebra of global sections of the {\em structure sheaf} $\SA=\SA_\CG$ on $\CG$. More generally, the global sections of any sheaf naturally form a $\CZ$-module (cf.~ \ref{Sheaves on moment graphs}). In this article we study the $\CZ$-modules that occur as global sections of sheaves. For this purpose we define a {\em localization functor} that is left adjoint to the functor of global sections. This gives rise to the notions of sheaves {\em generated by global sections} and of $\CZ$-modules {\em determined by local relations}. It turns out that the corresponding categories are equivalent (cf.~ \ref{combinatorial objects}). Hence we can view such objects either as sheaves on $\CG$ or as modules over $\CZ$.

The partial order defines the {\em Alexandrov topology} on $\CG$, i.e.~ the open sets are the sets that are downwardly closed. We consider the category $\CC^{fl}=\CC^{fl}(\CG)$ of sheaves that are generated by global sections and that are {\em flabby} for the Alexandrov topology (cf.~ Section \ref{section - categories of Z-modules}). They have the important property that their global sections can be constructed vertex by vertex following the order.

The structure sheaf  may fail to be flabby. As a substitute  Braden and MacPherson constructed the ``canonical sheaf'' $\SB=\SB_\CG$ on $\CG$ (it is called $\CM$ in \cite{BMP}). If $\CG$ comes from an equivariantly formal variety, $\SB$ encodes the structure of the equivariant and  ordinary intersection cohomologies in the sense that the space of global sections of $\SB$ is isomorphic to the hypercohomology of the equivariant $IC$-complex as a module over the equivariant cohomology ring. Specialization gives the analogous statement for the ordinary cohomology. Hence the failure of flabbiness of $\SA_\CG$ corresponds to the occurence of singularities. 

If $\CG$ is quasi-finite (cf.~ \ref{quasi-finite moment graphs}), any $\CZ$-module carries a functorial filtration indexed by the ordered set of vertices of $\CG$. If the subquotients of this filtration are free over $S$, we say that the module {\em admits a Verma flag}. If $\CG$ is a GKM-graph (cf.~ \ref{GKM graphs}), the corresponding sheaf is flabby. We denote by $\CV=\CC^{Verma}\subset\CC^{fl}$ the resulting category.  

Despite their universal nature it is not clear whether the Braden--MacPherson sheaves admit a Verma flag (cf.~ Remark \ref{proj and Verma flag}).  In the case that $\CG$ is associated to an equivariantly formal variety this is shown in \cite{BMP}. Our results confirm it if $\CG$ is associated to a non-critical block of category $\CO$. In general the subquotients of $\SB$ will only be reflexive. Let $\CC^{ref}\subset\CC^{fl}$ be the subcategory of such sheaves.

We show that the functorial filtration induces an exact structure on $\CV$ and on $\CC^{ref}$ and that the Braden--MacPherson sheaves associated to certain subgraphs are representatives of the indecomposable projective objects.

Let $\CO_\Lambda$ be a non-critical block of the {\em equivariant} category $\CO$ over a symmetrizable Kac--Moody algebra (cf.~ Section \ref{representation theory}), let $\CG_\Lambda$ be the associated moment graph and $\CV_\Lambda=\CV(\CG_\Lambda)$. Let $\CM_\Lambda\subset\CO_\Lambda$ be the subcategory of modules that admit a finite Verma flag. Our main result is that the structure functor (cf.~ \ref{structure functor}) gives rise to an equivalence 
$$
\CM_\Lambda\cong \CV_\Lambda
$$
of exact categories (Theorem \ref{equivalence of categories}). Under this equivalence the indecomposable projective modules correspond to the Braden-MacPherson sheaves. 
\subsection*{Acknowledgements}  This work was done during a stay at the University of Chicago that was financed by the Deutsche Forschungsgemeinschaft. I am grateful towards both institutions. I would like to thank  Tom Braden for helpful remarks on the first version of this paper.

\section{Sheaves on moment graphs}

In this section we define the basic notions such as moment graphs, sheaves on moment graphs, the structure sheaf and the functor of sections over a subgraph.
\subsection{Moment graphs}
Let $k$ be a field and $V$ a finite dimensional vector space over $k$. 
A {\em moment graph $\CG=(\CV,\CE,\leq,l)$ over V} is given by  
\begin{itemize}
\item a (not necessarily finite) graph $(\CV,\CE)$ with a set of vertices $\CV$ and a set of edges $\CE$,
\item a partial order $\leq$ on $\CV$ such that $x$ and $y$ are comparable if they are linked by an edge and
\item a labelling $l\colon\CE\to \DP^1(V)$ by one-dimensional subspaces of $V$.
\end{itemize}
We can think of the order as giving each edge a direction. We write $E\colon x\to y$ for an edge $E$ with endpoints $x$ and $y$ such that $x\leq y$. We write $E\colon x\xrightarrow{\alpha} y$, where $\alpha=\alpha(E)\in V^×$ is a generator of the line $l(E)$. If we forget about the direction we use the notations $E\colon x\linie  y$ and $E\colon x\stackrel{\alpha}\llinie y$. Often we write $x\in\CG$ for the vertices $x\in\CV$ of the underlying graph.

Later we will make two additional assumptions on the moment graph, namely quasi-finiteness (cf.~ \ref{quasi-finite moment graphs}) and the GKM-property (cf.~ \ref{GKM graphs}). For the Braden-MacPherson construction (cf.~ \ref{Braden-MacPherson sheaves}) we will in addition assume that the set of vertices is bounded from above.

\subsection{Examples of moment graphs}\label{examples of moment graphs}
\begin{enumerate}
\item The {\em generic moment graph} is the graph with one vertex and no edges. 
The {\em subgeneric moment graph} is the graph with two vertices that are connected by one labelled edge.
\item Let $(\CW,\CS)$ be a (not necessarily finite) Coxeter system and let $\fh$ be its real reflection representation (cf.~ \cite{Hum}). Let $\CT\subset\CW$ be the set of reflections in $\CW$. Let $\lambda\in\fh$ such that $\CW^\lambda:=\Stab_\CW(\lambda)$ is generated by simple reflections. Let $\CV=\CW.\lambda=\CW/\CW^\lambda$ be the orbit of $\lambda$. For any reflection $t\in\CT$ and $w.\lambda\in\CV$ with $tw.\lambda\neq w.\lambda$ define a labelled edge $E\colon w.\lambda\stackrel{w.\lambda-tw.\lambda}\lllinie tw.\lambda$. Together with the Bruhat order on $\CW/\CW^\lambda$ this defines the {\em $\fh$-moment graph $\CG(\CW,\CS,\CW^\lambda)$ associated to the Coxeter datum $(\CW,\CS,\CW^\lambda)$}.

\item To any block of the category $\CO$ over a symmetrizable complex Kac-Moody algebra outside the critical hyperplanes one can associate a moment graph over its Cartan subalgebra (cf.~ \ref{associated moment graph}). It is equivalent to the (complexified) moment graph associated to its Coxeter datum.

\item To a complex projective variety with a torus action that is equivariantly formal one associates a moment graph over the {\em dual} of the Lie algebra of the torus (cf.~ \cite{BMP}). It encodes the structure of the orbits of dimension zero and one. In the case of a flag variety of a complex finite dimensional Lie algebra we get the (complexified) moment graph associated to its Coxeter datum.
\end{enumerate}

\subsection{Sheaves on moment graphs}\label{Sheaves on moment graphs}
Let $\CG$ be a $V$-moment graph and $S=S(V)$ the symmetric algebra over $V$. We provide $S$ with an algebra $\DZ$-grading such that $\deg V=2$ (here this is completely arbitrary and chosen for reasons of compatibility with, for example, \cite{Soe3}). In the following we will consider graded $S$-modules and graded morphisms, i.e.~ morphisms of degree zero.

A {\em sheaf ${\SM}=(\{\SM^x\}, \{\SM^E\}, \{\rho_{x,E}\})$ on $\CG$} (or a {\em $\CG$-sheaf}\/) is given by 
\begin{itemize}
\item an $S$-module $\SM^x$ for any vertex $x\in\CV$,
\item an $S$-module $\SM^E$  for any edge $E\in\CE$ such that $l(E)\cdot \SM^E=0$,  
\item and maps $\rho_{x,E}\colon \SM^x\to \SM^E$ and $\rho_{y,E}\colon \SM^y\to\SM^E$ of $S$-modules if $x$ and $y$ are the  endpoints of $E$.  
\end{itemize}

A morphism $f\colon \SM\to \SN$ of $\CG$-sheaves consists of maps $f^x\colon \SM^x\to \SN^x$ and $f^E\colon \SM^E\to \SN^E$ that are compatible with $\rho$. The notion of a $\CG$-sheaf does not depend on the order on $\CV$. For a $\CG$-sheaf $\SM$  let $\supp(\SM):=\{x\in\CV\mid \SM^x\neq 0\}$ be its {\em support}. Let $\Sh(\CG)$ be the category of $\CG$-sheaves $\SM$ such that $\supp(\SM)$ is finite and each $\SM^x$ is finitely generated and torsion free over $S$.

Let $\SA=\SA_\CG$ be the {\em structure sheaf} of $\CG$, i.e.~ the sheaf defined by 
\begin{itemize}
\item
$\SA^x=S$ for a vertex $x$ of $\CG$,
\item 
$\SA^E=S/\alpha S$ for an edge $E\colon x\stackrel{\alpha}\llinie y$ and
\item $\rho_{x,E}$ the natural quotient map $S\to S/\alpha S$.
\end{itemize}

\subsection{Sections of sheaves on moment graphs}

A subgraph $\CI\subset \CG$ will always be  a full subgraph, i.e.~ any edge of $\CG$ connecting two vertices of $\CI$ will be contained in $\CI$. In other words, a subgraph is determined by its set of vertices. 

Let $\CI\subset\CG$ be a subgraph and $\SM$  a $\CG$-sheaf. The {\em space of sections of $\SM$ over $\CI$} is
$$
\SM(\CI):=
\left\{ 
(m_x)\in\prod_{x\in \CI}  \SM^x\left| \,
\begin{matrix}
\rho_{x,E}(m_x)=\rho_{y,E}(m_{y}) \\
\text{for any edge $E\colon x\linie y$}
\end{matrix}
\right.
\right\}.
$$
For subgraphs $\CJ\subset\CI\subset\CG$ we have a natural  restriction map
$\SM(\CI)\to \SM(\CJ)$. We write $\Gamma(\SM):=\SM(\CG)$ for the global sections of $\SM$.

Let 
$$
\CZ(\CI):=\SA(\CI)=
\left\{(z_x)\in\prod_{x\in\CI} S
\left|\, 
\begin{matrix}
z_x\equiv z_{y}\mod \alpha \\
\text{for any edge $x\stackrel{\alpha}\llinie y$ of $\CI$}
\end{matrix}
\right.
\right\}
$$
be the space of sections over $\CI$ of the structure sheaf of $\CG$. Pointwise multiplication defines an algebra structure on $\CZ(\CI)$. We call the algebra of global sections, $\CZ=\CZ(\CG)$, the {\em structure algebra} of $\CG$. 

\begin{remark} Suppose $\CG$ is the moment graph associated to a block of the category $\CO$ over a symmetrizable Kac-Moody algebra (cf.~ \ref{associated moment graph}). Let $R=S_{(\fh)}$ be the localization of $S=S(\fh)$ at the maximal ideal $S\cdot \fh$. Then the analogously defined algebra $\CZ_R$ is naturally isomorphic to the center of the block (cf.~ Theorem \ref{calculation of the center}). If $T$ is a torus and if $\CG$ is associated to an equivariantly formal $T$-variety, then $\CZ$ is naturally isomorphic to its $T$-equivariant cohomology (cf.~ \cite{GKM}).\end{remark}

For any $\CI\subset \CG$ let $\CZ^{\CI}$ be the image of the restriction $\CZ\to\CZ(\CI)$. Let $\CZ\catmod^f\subset\CZ\catmod$ be the full subcategory of $\CZ$-modules that are finitely generated over $S$, torsion free over $S$ and for which there is a {\em finite} subgraph $\CI\subset\CZ$ such that the action of $\CZ$ factors through $\CZ^\CI$.

If $\SM$ is a $\CG$-sheaf, then the space of its global sections $\Gamma(\SM)$ is a $\CZ$-module  under pointwise multiplication. If $\SM\in\Sh(\CG)$, then $\Gamma(\SM)\in\CZ\catmod^f$, hence we get a functor
$$
\Gamma\colon \Sh(\CG)\to \CZ\catmod^f.
$$

\subsection{Base change} \label{base change} Instead of the graded symmetric algebra $S=S(V)$ of the vector space we can also work with one of its localizations. Let $\fp\subset S$ be a graded ideal and let $S_\fp$ be the localization of $S$ at $\fp$. For any subgraph $\CI\subset\CG$ we define 
$$
\CZ_\fp(\CI):=\left\{ 
(z_x)\in\prod_{x\in \CI}  S_\fp\left| \,
\begin{matrix}
z_x\equiv z_{y} \mod \alpha\\
\text{for any edge $x\stackrel{\alpha}\llinie y$ of $\CI$}
\end{matrix}
\right.
\right\},
$$  
its {\em structure algebra over the base $S_\fp$}. Note that we have to consider only the edges $E$ such that $\alpha(E)\in S_\fp$ is not invertible. Hence it suffices to consider the {\em $\fp$-reduction} $\CI_\fp$ of $\CI$ that we construct from $\CI$ by keeping all vertices, but where we delete all edges labelled by elements that are invertible in $S_\fp$.

The canonical map $S\to S_\fp$ induces an algebra morphism $\CZ(\CI)\otimes_S S_\fp\to\CZ_\fp(\CI)$. 
If $\fp\subset S$ is such that $\CI_\fp=\CI$, i.e.~ such that $\alpha(E)\in \fp$ for all edges $E$, and if $\CI$ is finite, then we have a canonical isomorphism
$$
\CZ_\fp(\CI)=\CZ(\CI)\otimes_S S_\fp.
$$ 
If $\fp=S\cdot\gamma$ for some $\gamma\in V$ we will write $S_\gamma$ for $S_\fp$, $\CZ_\gamma(\CI)$ for $\CZ_\fp(\CI)$ and $\CI_\gamma$ for $\CI_\fp$. 

\section{Localization of $\CZ$-modules}
Let $\CG$ be a moment graph over $V$ and $\CZ=\CZ(\CG)$ its structure algebra. We want to define a functor $\Loc\colon \CZ\catmod^f\to\Sh(\CG)$ that is left adjoint to the functor $\Gamma$ of global sections. We call it the {\em localization functor}.
\subsection{Quasi-finite moment graphs}\label{quasi-finite moment graphs}
Let $\gamma\in V$. Recall that $S_\gamma$ is the localization of $S$ at the prime ideal $\fp=S\cdot\gamma$. Let $\CI\subset\CG$ be a finite subgraph.
Consider the composition $\CZ(\CG)^\CI\inj\CZ(\CI)\to\CZ_{\gamma}(\CI)$.
We call $\CG$ {\em quasi-finite} if the induced map
$$
\CZ(\CG)^\CI\otimes_S S_\gamma\to\CZ_{\gamma}(\CI)
$$
is a bijection for any {\em finite} subgraph $\CI\subset\CG$ and any $\gamma\in V$. Note that this map is always injective because it is an inclusion of subspaces of $\bigoplus_{x\in\CI}S_\gamma$.

\begin{lemma}\label{quasi-finiteness of graphs}
\begin{enumerate}
\item Any finite moment graph is quasi-finite.
\item Any moment graph associated to a Coxeter datum $(\CW,\CS,\CW^\lambda)$ (cf.~ \ref{examples of moment graphs}) is quasi-finite. 
\end{enumerate}
\end{lemma}
\begin{proof}

Let $\CG$ be as in $(1)$ or $(2)$, $V$ the underlying $k$-vector space, $\gamma\in V$ and choose a finite subgraph $\CI\subset\CG$. Then $\CZ_\gamma(\CI)=\bigoplus_\CJ\CZ_\gamma(\CJ)$, where $\CI_\gamma=\bigcup\CJ$ is the decomposition into connected components. Recall that we construct $\CI_\gamma$ from $\CI$ by deleting all edges that are not labelled by $k\cdot\gamma$. For any component $\CJ$ its structure algebra $\CZ_\gamma(\CJ)$ is generated, as an $S_\gamma$-module, by $d=(d_v)_{v\in\CI}$ with $d_v=1$ if $v\in\CJ$ and $d_v=0$ otherwise, and the set  $\{e^{w}\}_{w\in\CJ}$, where $e^w=(e^w_v)_{v\in\CI}$ with $e^w_w=\gamma$ and $e^w_v=0$ if $v\neq w$. Hence we have to show that  $d$ and $e^w$ are in the image of $\CZ(\CG)\otimes_S S_\gamma\to\CZ_\gamma(\CI)$.

Suppose $\CG$ is finite and let $\alpha_1,\dots,\alpha_n$ be the labels of all edges of $\CG$ that do not belong to $\CJ$ but that are connected to $\CJ$. In particular, none of the $\alpha_i$ is a multiple of $\gamma$. Define $a=(a_v)\in\bigoplus_{v\in\CG} S$ by $a_v=0$ if $v\not\in\CJ$ and $a_v=\alpha_1\cdots \alpha_n$ if $v\in\CJ$. Then $a\in\CZ(\CG)$ and $\alpha_1\cdots \alpha_n$ is invertible in $S_\gamma$, hence $(\alpha_1\cdots \alpha_n)^{-1}a\in\CZ(\CG)\otimes_S S_\gamma$ is a preimage of $d$.

Choose $w\in\CJ$. We can assume that there are at least two vertices in $\CJ$. Let $\alpha_1,\dots,\alpha_n$ be the labels of all edges connected to $w$. Define $b=(b_v)\in\bigoplus_{v\in\CG} S$ by $b_v=0$ for $v\neq w$ and $b_{w}=\gamma\cdot \beta_1\cdots \beta_m$, where $\beta_1,\dots,\beta_m$ are those of the $\alpha_i$ that are not a multiple of $\gamma$. Thus  $\beta_1\cdots \beta_m\in S_{\gamma}$ is invertible and $(\beta_1\cdots \beta_m)^{-1}b\in\CZ(\CG)\otimes_S S_\gamma$ is a preimage of $e^w$. Hence we proved the lemma for finite moment graphs.

Suppose $\CG=\CG(\CW,\CS,\CW^\lambda)$. Then every connected component of the reduced graph $\CI_\gamma$ is either generic or subgeneric, since there is at most one edge at any vertex labelled by a multiple of $\gamma$. 
Recall that $V=\fh$ is the reflection representation of $\CW$.
The element $c^\lambda=(c^\lambda_{\bar v})\in\prod_{\bar v\in\CW/\CW^\lambda}S$ with $c^\lambda_{\bar v}:=\bar v(\lambda)\in\fh$
is in $\CZ(\CG)$, and $c^\lambda-\bar v(\lambda)\cdot 1$ is supported on $\CG\setminus\{\bar v\}$.

We assume that $\CJ=\bar w\stackrel{\gamma}\llinie\bar w^\prime$ is subgeneric (the generic case is proven similarly).
 Then $\bar w^\prime=t\bar w$ for some reflection $t$ and $\gamma$ is a multiple of $\bar w(\lambda)-t\bar w(\lambda)$. We have to show that the elements $(\gamma,0)$ and $(1,1)$ of $\CZ_\gamma(\CJ)\subset\CZ_\gamma(\CI)$ have a preimage in $\CZ(\CG)\otimes_S S_\gamma$.

Consider $a=\prod_{\bar v\in\CI,\bar v\neq\bar w}(c^\lambda-\bar v(\lambda)\cdot 1)\in\CZ(\CG)$. Then $\supp(a)\cap\CI=\{\bar w\}$. For $\bar v\not\in\{\bar w,t\bar w\}$ the element $c_{\bar w}^\lambda-\bar v(\lambda)=\bar w(\lambda)-\bar v(\lambda)$ is not divisible by $\gamma$, hence invertible in $S_\gamma$. Moreover, $c^\lambda_{\bar w}-t\bar w(\lambda)$ is a scalar multiple of $\gamma$. Thus multiplying with an invertible element gives a preimage in $\CZ(\CG)\otimes_S S_\gamma$ of $(\gamma,0)\in\CZ_\gamma(\CJ)$. 

Consider $b=\prod_{\bar v\in\CI,\bar v\not\in\{\bar w,t\bar w\}}(c^\lambda-\bar v(\lambda)\cdot 1)\in\CZ(\CG)$. Then $\supp(b)\cap \CI=\{\bar w,t\bar w\}$. Neither $b_{\bar w}$ nor $b_{t\bar w}$ is divisible by $\gamma$, but $b_{\bar w}-b_{t\bar w}$ is divisible by $\gamma$. We already showed that $(\gamma,0)$ lies in the image of $\CZ(\CG)^\CI$, hence there is a preimage of  $(x,x)\in\CZ_\gamma(\CJ)$ for some invertible element $x\in S_\gamma$. Thus also $(1,1)$ has a preimage. 
\end{proof}

\subsection{Modules over the structure algebra}\label{decomposition at generic point}
Suppose that $\CG$ is quasi-finite.
Let $Q$ be the quotient field of $S$. 
For any $S$-module $N$ let $N_Q:=N\otimes_S Q$.  We have $\CZ(\CG)^\CI\otimes_S Q=\bigoplus_{x\in\CI}Q$ for any finite subgraph $\CI\subset\CG$. So for $M\in\CZ\catmod^f$ there is a canonical action of $\prod_{x\in\CG} Q$ on $M_Q$, hence a canonical decomposition
$$
M_Q=\bigoplus_{x\in\CG} M^x_{Q}.
$$
Since we assume that $M$ has no $S$-torsion we get an inclusion $M\hookrightarrow M_Q=\bigoplus_{x\in\CG} M^x_{Q}$, so we can view any $m\in M$ as a tuple $m=(m_x)$ with $m_x\in M^x_{Q}$. We define $\supp(m):=\{x\in\CG\mid m_x\neq 0\}$ and $\supp(M):=\bigcup_{m\in M} \supp(m)\subset\CG$. If the action of $\CZ$ on $M$ factors over $\CZ^\CI$, then $\supp(M)\subset\CI$, hence $\supp(M)$  is finite.

For any $\CJ\subset\CG$ let 
\begin{itemize}
\item $M_\CJ$ be the subspace of elements supported on $\CJ$, i.e.~
$M_\CJ=M\cap \bigoplus_{x\in\CJ} M_Q^x$,
\item $M^\CJ$  the projection of $M$ onto the coordinates in $\CJ$, i.e.~ let $M^\CJ$ be the image of $M$ under the projection $\bigoplus_{x\in\CG} M_Q^x\to \bigoplus_{x\in\CJ} M_Q^x$ with kernel $\bigoplus_{x\in\CG\setminus \CJ} M_Q^x$.
\end{itemize}
Note that $M^\CJ$ and $M_{\CJ}$ are $\CZ$-modules.
We write $M^{x_1,\dots, x_n}$ and $M_{x_1,\dots, x_n}$ instead of $M^{\{x_1,\dots, x_n\}}$ and $M_{\{x_1,\dots, x_n\}}$. For any $\CJ\subset\CG$ we have inclusions
$$
\bigoplus_{x\in\CJ} M_x\subset M_\CJ\subset 
M^{\CJ}\subset\bigoplus_{x\in\CJ} M^x.
$$

\subsection{Localization} For any $M\in\CZ\catmod^f$
we want to construct a $\CG$-sheaf $\Loc(M)=(\{\Loc(M)^x\},\{\Loc(M)^E\},\rho_{x,E})$.  For a vertex $x$ of $\CG$ we define
$$
\Loc(M)^x:=M^x.
$$ 
For any edge $E\colon x\stackrel{\alpha}\llinie y$  of $\CG$  let $\CZ(E):=\{(z_x,z_y)\in S\oplus S\mid z_x\equiv z_y\mod \alpha\}$ be its local structure algebra. Note that the restriction map $\CZ\to\CZ(E)$ is not necessarily surjective. Consider the inclusion $M^{x,y}\subset M^x\oplus M^y$. The algebra $\CZ(E)$ acts on $M^x\oplus M^y$ by pointwise multiplication. Let $M(E):=\CZ(E)\cdot M^{x,y}$  be the submodule generated by $M^{x,y}$ and let $M(E)\to M^x$ and $M(E)\to M^y$ be the projections onto the coordinates. Let $\Loc(M)^E$ be the push out in the following diagram  of $S$-modules 

\centerline{
\xymatrix{
M(E)\ar[d]\ar[r] &  M^x\ar[d] \\
M^y\ar[r]&\Loc(M)^E 
}
}
\noindent
and let $\rho_{x,E}\colon \Loc(M)^x=M^x\to \Loc(M)^E$ and $\rho_{y,E}\colon \Loc(M)^y=M^y\to \Loc(M)^E$ be the push out maps. Since $M(E)\to M^x$ and $M(E)\to M^y$ are surjective they induce isomorphisms 
$$
M^x/M(E)_x\xrightarrow{\sim} \Loc(M)^E\xleftarrow{\sim} M^y/M(E)_y,
$$
where $M(E)_x=\ker(M(E)\to M^y)\cong M^x\cap M(E)$ and $M(E)_y=\ker(M(E)\to M^x)\cong M^y\cap M(E)$.

Since $(\alpha,\alpha)=(0,\alpha)+(\alpha,0)\in\CZ(E)$ we have $\alpha M^x=(\alpha,0)\cdot M(E)\subset M(E)_x$ and hence $\Loc(M)^E$ is annihilated by $\alpha$. So we have defined a $\CG$-sheaf $\Loc(M)$. 
 Its support is finite and the $S$-modules on the vertices are torsion-free and finitely generated, hence $\Loc(M)\in\Sh(\CG)$. This gives the {\em localization functor}
 $$
\Loc\colon \CZ\catmod^f\to \Sh(\CG).
$$
\subsection{Properties of the localization functor} Let $M\in\CZ\catmod^f$, $E\colon x\linie y$ an edge and consider the local sections
$$
\Loc(M)(E)=\{(m_x,m_y)\in M^x\oplus M^y\mid \rho_{x,E}(m_x)=\rho_{y,E}(m_y)\}.
$$ 
For later use we state the following two lemmas.
\begin{lemma}\label{subgeneric localization} We have
$M(E)=\Loc(M)(E)$
as subspaces of $M^x\oplus M^y$.
Hence globally we have
$$
\Gamma\circ\Loc(M)=\left\{ (m_x)\in\bigoplus_{x\in\CG} M^x\left| 
\,
\begin{matrix}
(m_x,m_y)\in M(E) \\ 
\text{\rm for any edge $E\colon x\linie y$}
\end{matrix}
\right.
\right\}.
$$
\end{lemma}

\begin{proof} We obviously have $M(E)\subset\Loc(M)(E)$. Let $(m_x,m_y)\in\Loc(M)(E)$. We want to show that $(m_x,m_y)\in M(E)$. There is $m^\prime=(m^\prime_x,m^\prime_y)\in M(E)$ such that $m^\prime_x=m_x$, so we can assume that $m_x=0$. Then the image of $m_y$ in $\Loc(M)^E$ is zero. Since $\Loc(M)^E\cong M^y/M(E)_y$ we deduce $(0,m_y)\in M(E)_y\subset M(E)$.   
\end{proof}

\begin{lemma}\label{sections and restriction} 
For any $M\in\CZ\catmod^f$ and any $\CI\subset \CG$ we have $
\Loc(M)(\CI)= \Loc(M^\CI)(\CG)$.
\end{lemma}
\begin{proof} For any $\CG$-sheaf $\SN$ let $\SN^\CI$ be its restriction to the subgraph $\CI$, i.e.~ the $\CG$-sheaf that we obtain from $\SN$ by setting all $\SN^x$ and $\SN^E$ to zero if $x$ or $E$ do not belong to $\CI$. Then
$\Loc(M)(\CI)  =  \Loc(M)^\CI(\CG)$.
By construction we have $\Loc(M)^\CI=\Loc(M^\CI)$, hence the claim. 
\end{proof}

\subsection{Natural transformations} 
We want to define natural transformations $f\colon \Loc\circ\Gamma\to\id_{\Sh(\CG)}$ 
and $g\colon \id_{\CZ\catmod^f}\to \Gamma\circ \Loc$.

 Let $\SM\in\Sh(\CG)$. We have to define the components $f^x\colon \Loc(\Gamma(\SM))^x\to \SM^x$ and $f^E\colon \Loc(\Gamma(\SM))^E\to \SM^E$. The space of global sections $\Gamma(\SM)$ is a subspace of $\bigoplus_{x\in\CG} \SM^x$ and $\CZ$ acts pointwise, hence we have an inclusion $f^x\colon\Loc(\Gamma(\SM))=\Gamma(\SM)^x\subset \SM^x$. 

Let $E\colon x\linie y$ be an edge. The diagram

\centerline{
\xymatrix{
\Gamma(\SM)(E) \ar[r]\ar[d]& \Gamma(\SM)^x \ar[d]^{\rho_{x,E}} \\
\Gamma(\SM)^y \ar[r]^{\rho_{y,E}} & \SM^E
}}
\noindent
commutes, so the universal property gives a map $f^E\colon\Loc(\Gamma(\SM))^E\to \SM^E$.
It is clear that $f=(f^x, f^E)$ is compatible with $\rho$, so we have defined a natural transformation $f\colon \Loc\circ \Gamma\to \id_{\Sh(\CG)}$. 

Let $M\in\CZ\catmod^f$. Then $\Loc(M)^x=M^x$ and any $m=(m_x)\in M$ defines a global section in $\Gamma(\Loc(M))$ by Lemma \ref{subgeneric localization}, hence we get an inclusion $M\to \Gamma(\Loc(M))$ and a natural transformation
$g\colon \id_{\CZ\catmod^f}\to \Gamma\circ \Loc$.

\subsection{Adjointness} Let $\CA$, $\CB$ and $\CC$ be categories, $F\colon \CA\to \CB$ and $G,H\colon \CB\to \CC$ functors and $j\colon G\to H$ a natural transformation. Then, composing with $F$, we get an induced transformation between the functors $G\circ F,H\circ F\colon \CA\to\CC$ that we denote by $j\diamond F$. In a similar situation we define $F\diamond j$ analogously.

\begin{proposition}\label{Isomorphisms of functors} The natural transformations
$\Gamma\xrightarrow{g\diamond\Gamma} \Gamma\circ \Loc \circ\Gamma$
and $\Loc\circ\Gamma\circ \Loc \xrightarrow{f\diamond \Loc} \Loc$
are isomorphisms with inverses $\Gamma\diamond f$ and $\Loc\diamond g$, resp.
\end{proposition}
\begin{proof} Let $\SM\in\Sh(\CG)$. Then 
$$
\Gamma(\SM)\xrightarrow{g\diamond\Gamma}\Gamma\circ\Loc\circ\Gamma(\SM)\xrightarrow{\Gamma\diamond f}\Gamma(\SM)
$$
are mutual inclusions of subspaces of $\bigoplus_{x\in\CG}\SM^x$, hence the first claim.

Let $M\in\CZ\catmod^f$ and $x\in\CG$. Then $(\Gamma\circ \Loc)(M)^x\subset M^x$ and composing with $M^x\xrightarrow{g^x} (\Gamma\circ\Loc)(M)^x$ gives the identity.  Applying $\Loc$ yields
$$
\Loc(M)^x\xrightarrow{\Loc\diamond g}(\Loc\circ\Gamma\circ\Loc)(M)^x\xrightarrow{f\diamond\Loc}\Loc(M)^x.
$$  

Let $E\colon x\linie y$ be an edge. Then $M^{x,y}\xrightarrow{g^{x,y}}(\Gamma\circ \Loc)(M)^{x,y}$ is an inclusion of subspaces of $M^x\oplus M^y$. Every global section specializes to a local section, hence $(\Gamma\circ \Loc)(M)^{x,y}\subset\Loc(M)(E)$. Acting with $\CZ(E)$ gives
$$
M(E)\subset(\Gamma\circ \Loc)(M)(E)\subset\Loc(M)(E).
$$
By Lemma \ref{subgeneric localization} these are bijections, hence $M(E)=(\Gamma\circ\Loc)(M)(E)$. Then $\Loc(M)^E=(\Loc\circ\Gamma\circ\Loc)(M)^E$ and $(\Loc\diamond g)^E$ and $(f\diamond \Loc)^E$ are mutually inverse isomorphisms.
\end{proof}

\begin{theorem}
The functor $\Loc\colon\CZ(\CG)\catmod^f\to\Sh(\CG)$ is left adjoint to the functor $\Gamma\colon\Sh(\CG)\to\CZ(\CG)\catmod^f$.
\end{theorem}
\begin{proof}
The maps $f$ and $g$ induce transformations $\Hom(\cdot,\Gamma\cdot)\to\Hom(\Loc\cdot,\cdot)$ and $\Hom(\Loc\cdot,\cdot)\to\Hom(\cdot,\Gamma\cdot)$
of bifunctors, respectively. In order to show that they are mutually inverse, we have to show that the two compositions $\Loc\xrightarrow{\Loc\diamond g} \Loc\circ\Gamma\circ \Loc \xrightarrow{f\diamond \Loc} \Loc$
and $
\Gamma\xrightarrow{g\diamond\Gamma} \Gamma\circ \Loc \circ\Gamma\xrightarrow{\Gamma\diamond f} \Gamma$
both equal the identity. The statement of Proposition \ref{Isomorphisms of functors} is even stronger.
\end{proof}

\subsection{Combinatorial objects}\label{combinatorial objects}

From the Proposition \ref{Isomorphisms of functors} it follows that $\Gamma$ and $\Loc$ induce isomorphisms between their respective essential images.  Concretely, this means the following. Let $\CZ\catmod^{loc}$ be the full subcategory of $\CZ\catmod^f$ consisting of all objects that are ``determined by local relations", i.e.~ the category of objects $M\in\CZ\catmod^f$ such that the natural transformation $(\Gamma\circ \Loc)(M)\to M$ is an isomorphism. Analogously, let $\Sh(\CG)^{glob}$ be the full subcategory of $\Sh(\CG)$ consisting of all objects that are ``generated by global sections", i.e.~ the category of objects $\SM\in\Sh(\CG)$ such that the natural transformation $\SM\to (\Loc\circ\Gamma)(\SM)$ is an isomorphism. By Proposition \ref{Isomorphisms of functors} the image of $\Loc$ lies in $\Sh(\CG)^{glob}$, the image of $\Gamma$ lies in $\CZ\catmod^{loc}$ and $\Loc$ and $\Gamma$ define mutually inverse equivalences between $\CZ\catmod^{loc}$ and $\Sh(\CG)^{glob}$. We can think of this as a category $\CC=\CC(\CG)$ that is simultaneously embedded into $\Sh(\CG)$ and $\CZ\catmod^f$, so each object of $\CC$ has both a local and a global nature. 

\section{Categories of $\CZ$-modules and exact sequences}\label{section - categories of Z-modules}
Let $\CG$ be a quasi-finite moment graph and $\CZ=\CZ(\CG)$ its structure algebra. 
We can topologize the set of vertices of $\CG$ using the Alexandrov topology on partially ordered sets. That means that the open sets are the subsets
$\CI\subset\CG$ that are downwardly closed, i.e.~ if $y\in\CI$ and $x\leq y$, then $x\in\CI$.

Let $M\in\CZ\catmod^{f}$. We say that $M$ is {\em flabby}, if $M^\CI$ is determined by local relations, i.e.~ $M^\CI\in\CZ\catmod^{loc}$, for any open subgraph $\CI\subset\CG$. In particular, $M=M^\CG$ is determined by local relations. Proposition \ref{flabby sheaves} shows that the associated sheaf is flabby for the Alexandrov topology. Let
$$
\CZ\catmod^{fl}\subset\CZ\catmod^{loc}
$$
be the full subcategory of flabby objects.

Let $M=\bigoplus_{n\in\DZ} M_n$ be an $S$-module, graded as always. For $i\in\DZ$ define the graded $S$-module $M[i]$ by $M[i]_n=M_{n+i}$. For graded $S$-modules $M$ and $N$ let $\Hom^i_S(M,N):=\Hom_S(M,N[i])$ and set $\Hom_S^\bullet(M,N)=\bigoplus_{i\in\DZ}\Hom^i_S(M,N)$. The graded dual of $M$ is the graded $S$-module $M^\ast=\Hom^\bullet_S(M,S)$. Note that $S^\ast=S$ canonically. Recall that $M$ is called {\em reflexive}, if it is isomorphic to its double dual, i.e.~ if $M\xrightarrow{\sim} (M^\ast)^\ast$. This is the case if and only if $M$ is the intersection of all its graded localizations at graded prime ideals of $S$ of height one. 

Let $M\in\CZ\catmod^f$. We say that $M$ is  {\em reflexive}, if $M^\CI$ is reflexive as an $S$-module  for any open subgraph $\CI\subset\CG$. Let
$$
\CZ\catmod^{ref}\subset\CZ\catmod^f
$$
be the full subcategory of reflexive objects. 

We say that $M$ {\em admits a Verma flag},  if $M^\CI$ is {\em graded free}, i.e.~ isomorphic to $\bigoplus S[l_i]$ for some $l_1,\dots,l_n\in\DZ$, for any open subgraph $\CI\subset\CG$. Let
$$
\CZ\catmod^{Verma}\subset\CZ\catmod^{ref}
$$
be the full subcategory of objects that admit a Verma flag.

In the following we will define an exact structure (i.e.~ a class of exact sequences) on $\CZ\catmod^f$. The subcategories $\CZ\catmod^{ref}$ and $\CZ\catmod^{Verma}$ are closed under extensions in $\CZ\catmod^f$, hence inherit an exact structure. Since $\CZ\catmod^{fl}\subset\CZ\catmod^f$ is not closed under extensions it is not clear if it inherits an exact structure. I do not know if it does. 

If $\CG$ is  a GKM-graph (see \ref{GKM graphs} for the definition) we will show that any reflexive $\CZ$-module is flabby, i.e.~ that
$$
\CZ\catmod^{Verma}\subset\CZ\catmod^{ref}\subset\CZ\catmod^{fl}\subset\CZ\catmod^{loc}.
$$
Thus we can consider all of the categories above as categories of sheaves on $\CG$. To incorporate this viewpoint we denote them by 
$$
\CC^{Verma}\subset\CC^{ref}\subset\CC^{fl}\subset\CC.
$$

The most important category for us is the category $\CV:=\CC^{Verma}$ of objects admitting a Verma flag. In Section \ref{representation theory} we prove the main result of this article, namely that if $\CG$ is associated to a non-critical block of the equivariant category $\CO$ over a symmetrizable Kac-Moody algebra, then $\CV$ is equivalent (as an exact category) to the subcategory of modules that admit a Verma flag. 

\subsection{An exact structure}\label{subsection - exact structure}

Let $A \to B\to C$
be a sequence in $\CZ\catmod^f$. We say that it is {\em short exact}, if for any open $\CI\subset\CG$ the induced sequence 
$$
0\to A^{\CI} \to B^{\CI}\to C^{\CI} \to 0
$$
is a short exact sequence of abelian groups.

In particular any short exact sequence $A\to B\to C$ induces a short exact sequence $0\to A\to B\to C\to 0$ of $\CZ$-modules. Conversely, if $0\to A\to B\to C\to 0$ is short exact as a sequence of $\CZ$-modules, then $A^\CI\to B^\CI$ is injective, $B^\CI\to C^\CI$ is surjective and  we have to check the exactness of $0\to A^\CI\to B^\CI\to C^\CI\to 0$  only in the middle.

\begin{theorem}\label{Theorem - exact structure} The above definition gives an exact structure in the sense of Quillen (cf.~ \cite{Qu}) on  $\CZ\catmod^f$. 
\end{theorem}

\begin{proof}
Let $E$ be the class of short exact sequences in $\CZ\catmod^f$. 
A morphism $g\colon B\to C$ in $\CZ\catmod^f$ is called an admissible epimorphism if it occurs on the right side of a short exact sequence. Dually, a morphism $f\colon A\to B$ is called an admissible monomorphism if it occurs on the left side of an exact sequence. We have to show the following properties of a Quillen exact structure:

\begin{enumerate}
\item Any sequence isomorphic to a sequence in $E$ is in $E$. Any split sequence is in $E$. If $A\xrightarrow{f} B\xrightarrow{g} C$ is in $E$, then $f$ is a kernel of $g$ and $g$ is a cokernel of $f$ in $\CZ\catmod^f$. 
\item The class of admissible epimorphisms is stable under composition and under base-change along arbitrary morphisms in $\CZ\catmod^f$. Dually, the class of admissible monomorphisms is stable under composition and under cobase-change along arbitrary morphisms in $\CZ\catmod^f$.
\item If $g\colon B\to C$ has a kernel in $\CZ\catmod^f$ and if there is a morphism $h\colon D\to B$ such that $g\circ h\colon D\to C$ is an admissible epimorphism, then $g$ is an admissible epimorphism. Dually for admissible monomorphisms. 
 \end{enumerate}
The properties $(1)$ follow immediately from the definition. We will now show that admissible monomorphisms are stable under composition. 

Let $A\to B\to C $ and $B\to D\to E$ be short exact sequences in $\CZ\catmod^f$. Let $F$ be the cokernel of $A\to D$ in the category of $\CZ$-modules. Then we have a short exact sequence (of $\CZ$-modules) $0\to C\to F\to E\to 0$, hence $F$ is finitely generated and torsion free over $S$, so $F\in\CZ\catmod^f$. We have to show that the sequence $A\to D\to F$ is short exact, i.e.~ we have to check if the sequence $0\to A^\CI\to D^\CI\to F^\CI\to 0$ is short exact in the middle for any open subgraph $\CI\subset \CG$. This follows from the diagram
$$
\begin{matrix}
A^{\CI} & \to & B^{\CI} & \to & C^{\CI}\\
\parallel && \downarrow &&\downarrow \\
A^{\CI} & \to & D^{\CI} & \to & F^{\CI} \\
&& \downarrow && \downarrow \\
&& E^{\CI} & = & E^{\CI} 
\end{matrix}
$$
using the fact that the middle row and the top column are exact and that the map $C^{\CI}\to F^\CI$ is injective. That the composition of two admissible epimorphisms is an epimorphism is proven in a similar way.  

Now we consider base-change. Let $A\to B\to C$ be a short exact sequence and $D\to C$ an arbitrary morphism. Then there is a base-change diagram 
$$
\begin{matrix}
A & \to & E & \to & D \\
\parallel & & \downarrow && \downarrow \\
A & \to & B & \to & C
\end{matrix}
$$
in $\CZ\catmod^f$.
We have to show that the top row is short exact. Again we only have to show that for any open $\CI\subset\CG$ the sequence $0\to A^\CI\to E^\CI \to D^\CI\to 0$ is exact in the middle. This follows from the diagram
$$
\begin{matrix}
&&A^\CI & \to & E^\CI & \to & D^\CI \\
&&\parallel & & \downarrow && \downarrow \\
0 & \to & A^\CI & \to & B^\CI & \to & C^\CI
\end{matrix}
$$
using the fact that the bottom row is exact and the injectivity of the map $E^\CI\to B^\CI\oplus D^\CI$. Cobase-change for admissible epimorphisms is treated similarly. 

According to Keller (cf.~ \cite{Kel}), (3) follows from (1) and (2). 
\end{proof}

\subsection{Flabby objects}
We want to characterize the flabby $\CZ$-modules when considered as sheaves on the moment graph. 
Let $x$ be a vertex of $\CG$ and let $\CE_{\delta x}\subset\CE$ be the set of all edges $E\colon y\to x$ ending at $x$. For a $\CG$-sheaf $\SM$ let $\SM^{\delta x}\subset\bigoplus_{E\in\CE_{\delta x}}\SM^E$ be the image of the map
$$ 
\SM(\{< x\})\subset \bigoplus_{y<x}\SM^y\xrightarrow{\pi}  \bigoplus_{y\to x} \SM^y \xrightarrow{\bigoplus \rho_{y,E}} \bigoplus_{E\in\CE_{\delta x}}\SM^E,
$$ 
where $\pi$ is the projection with the sum of the $\SM^y$ with $y<x$ and  $y$ not linked to $x$ as the kernel.
Note that for $\SM\in\Sh(\CG)^{glob}$ the image of $\SM^x\xrightarrow{\bigoplus \rho_{x,E}} \bigoplus_{E\in\CE_{\delta x}} \SM^E$
is contained in $\SM^{\delta x}$ since any $m_x\in \SM^x$ comes from a global section.

\begin{proposition} \label{flabby sheaves} 
Let $\SM\in\Sh(\CG)^{glob}$. Then the following are equivalent:
\begin{enumerate} 
\item $\Gamma(\SM)\in\CZ\catmod^{fl}$, i.e.~ $\Gamma(\SM)$ is flabby.
\item The sheaf $\SM$ is flabby for the Alexandrov topology on the set of vertices of $\CG$, i.e.~ for any open $\CI\subset\CG$ the restriction map
$\SM(\CG)\to \SM(\CI)$
is surjective.
\item For any vertex $x\in\CG$ the restriction map
$\SM(\{\leq x\})\to \SM(\{<x\})$
is surjective.
\item For any vertex $x$ of $\CG$ the map
$\bigoplus_{E\in\CE_{\delta x}} \rho_{x,E}\colon
\SM^x\to \SM^{\delta x}$
is surjective. 
\end{enumerate}
\end{proposition}

\begin{proof} Let $i_1\colon \Gamma(\SM)^\CI\inj (\Gamma\circ \Loc)(\Gamma(\SM)^\CI)$ be induced by the natural transformation $\id\to \Gamma\circ \Loc $, and $i_2\colon \Gamma(\SM)^\CI\inj \SM(\CI)$  the canonical inclusion. By Lemma \ref{sections and restriction} 
 we can identify $(\Gamma\circ \Loc)(\Gamma(\SM)^\CI)\cong \Loc\circ\Gamma(\SM)(\CI)= \SM(\CI)$ and, under these identifications, $i_1=i_2$. Property (1) means that $i_1$ is an isomorphism for any $\CI$, whereas property (2) says that $i_2$ is an isomorphism for any $\CI$. Hence (1) and (2) are equivalent.

It is clear that (3) and (4) are equivalent. And (3) follows from a special case of (2). It remains to prove that (3) implies (2).  

It is enough to show that, if (3) holds and if $\CI^\prime=\CI\dot\cup\{x\}$ (for some vertex $x$) is again open, then the restriction
$\SM(\CI^\prime)\to \SM(\CI)$ 
is surjective.  Let $m\in\SM(\CI)$ and let $m^{<x}\in\SM(\{{<x}\})$ be its restriction. By (3) we find an extension $m^\prime\in\SM(\{{\leq x}\})$ of $m^{<x}$. Then $m^{\prime\prime}=(m^{\prime\prime}_y)$ with $m^{\prime\prime}_y=m_y$ if $y\neq x$ and $m^{\prime\prime}_x=m^\prime_x$ is a preimage of $m$ in $\SM(\CI^\prime)$. 
\end{proof}

\subsection{The order filtration}
Let $M\in\CZ\catmod^{fl}$ and let $x$ be a vertex of $\CG$. We denote by $M^{\leq x}$ the module $M^\CI$ for $\CI=\{y\mid y\leq x\}$. Similarly we define $M^{<x}$. For $y\leq x$ we have a  map $M^{\leq x}\to M^{\leq y}$, so we get a filtration of $M$ indexed by the partially ordered set of vertices. We call it the {\em  order filtration}. For a vertex $x$ we denote by $M^{[x]}$ the kernel of this filtration at the index $x$, i.e.~ 
$$
M^{[x]}:=\ker(M^{\leq x}\to M^{<x}).
$$
Note that $M^{[x]}\subset M^{\leq x}$ is the submodule of elements supported on $\{x\}$, hence $M^{[x]}$ is the intersection of the kernels of the maps $M^{\leq x}\to M^{\leq y}$ for all $y<x$. If $\SM$ is the sheaf associated to $M$, then  $M^{\varle x}=\CM(\{\varle x\})$ because $M$ is flabby, hence $M^{[x]}=\SM^{[x]}:=\ker(\SM^x\to \SM^{\delta x})$. Let $\CI\subset\CG$ be open and $x\in\CI$ maximal. Since $x$ is not linked to any vertex in $\CI\setminus\{\leq x\}$ we get the following lemma.

\begin{lemma}\label{ext of Vermas}  
There is a canonical exact sequence 
$$
0\to M^{[x]}\to M^{\CI} \to M^{\CI\setminus\{x\}}\to 0.
$$
Thus $M$ is an extension of the $M^{[x]}$ ``in the direction given by the order'', i.e.~ for any enumeration  $i\colon\supp(M)\to\{1,2,\dots,\#\supp(M)\}$ such that $i(v)<i(w)$ if $w<v$ there is a filtration by $\CZ$-modules
$$
0=M_{0}\subset M_{1}\subset\dots\subset M_{\#\supp(M)}=M
$$
with $M_{i(v)}/M_{i(v)-1}\cong M^{[v]}$ for all $v\in\supp(M)$.
\end{lemma}

The following lemma is an immediate consequence.

\begin{lemma} \label{exact sequences} Let $A \to B\to C$ be a sequence in $\CZ\catmod^f$ with $A,B,C\in\CZ\catmod^{fl}$. It is short exact if and only if for any $x\in\CG$ the sequence 
$$
0\to A^{[x]}\to B^{[x]}\to C^{[x]}\to 0
$$
is a short exact sequence of $\CZ$-modules.
\end{lemma}

 \subsection{GKM-graphs and reflexive $\CZ$-modules}\label{GKM graphs}
A moment graph $\CG$ is called a  {\em Goresky-Kottwitz-MacPherson}-graph (or {\em GKM-graph}) if it is quasi-finite and  if, for any vertex $x$, the labels of the edges that are linked to $x$ are pairwise linearly independent, i.e.~  $\CG_\fp$ splits into a disjoint union of generic and subgeneric graphs for any prime $\fp\subset S$ of height one.  

We want to show that if $\CG$ is a GKM-graph, then every reflexive $\CZ$-module $M$ is determined by local relations. Applying this to $M^\CI$ for any open subgraph $\CI\subset\CG$ shows that $\CZ\catmod^{ref}$ is a subcategory of $\CZ\catmod^{fl}$. Being reflexive as well as being free is stable under extensions (in the category of $S$-modules), hence $\CZ\catmod^{ref}$ and $\CZ\catmod^{Verma}$ inherit an exact structure from $\CZ\catmod^f$.

From now on we assume that $\CG$ is a GKM-graph.  
If $M\in\CZ\catmod^f$, then $M_\fp$ can be considered as a $\CZ_{\fp}(\CG_\fp)$-module by quasi-finiteness. 
\begin{proposition} Suppose that $\CG$ is a GKM-graph and that $M\in\CZ\catmod^f$ is reflexive as an $S$-module. Then $M$ is determined by local relations, i.e.~ $M\in\CZ\catmod^{loc}$.\end{proposition}

\begin{proof}
We have to show that the map $M\to (\Gamma\circ \Loc)(M)$ is an isomorphism. Let $\fp\subset S$ be a prime ideal of height one.
Let $\Gamma_\fp$ and $\Loc_\fp$ be the global section and localization functors for $\CG_\fp$ (over the base $S_\fp$). Since $\CG$ is quasi-finite, $M_\fp$ decomposes into the direct sum of its submodules supported on the connected components of $\CG_\fp$. By assumption, $\CG_\fp$ is a disjoint union of generic and subgeneric moment graphs, and hence Lemma \ref{subgeneric localization} shows that the map $M_\fp \to (\Gamma_\fp\circ \Loc_\fp) (M_\fp)$ is an isomorphism. The natural transformation $\id\to \Gamma\circ \Loc $ is compatible with localizations. The following commutative diagram

\centerline{
\xymatrix{
{{M}} \ar[d]_{\wr} \ar @{^{(}->}[r] & {{(\Gamma\circ\Loc)(M)}} \ar @{^{(}->}[d] \\
\bigcap M_\fp\ar[r]^{\sim\hspace{30pt}} &  \bigcap(\Gamma_\fp\circ \Loc_\fp)(M_\fp)
}
}
\noindent
proves the proposition (the intersections are taken over all prime ideals of height one).
\end{proof}

\begin{lemma}\label{sections on reflexives} Let $M\in \CZ\catmod^{fl}$ be flabby and suppose that for any $x\in\CG$ the $S$-module $M^x$ is reflexive. Then $M\in\CZ\catmod^{ref}$. 
\end{lemma}
\begin{proof} Let $\CI\subset\CG$ be an open subset.  By Proposition \ref{flabby sheaves}, (2) we have $M^\CI=(\Gamma\circ\Loc)(M)^\CI=\Loc(M)(\CI)$. Then
\begin{eqnarray*}
\Loc(M)(\CI) & \subset &  \bigcap_{\height(\fp)=1} \Loc(M)(\CI)_\fp \\
& \subset &  \bigcap_{\height(\fp)=1} 
\left\{ 
\{m_x\}\in\bigoplus_{x\in \CI}  M_\fp^x\left| \,
\begin{matrix}
\rho_{x,E}(m_x)=\rho_{y,E}(m_{y}) \\
\text{for any edge $E\colon x\linie y$}
\end{matrix}
\right.
\right\} \\
& = & \Loc(M)(\CI).
\end{eqnarray*}
Thus $M^\CI=\Loc(M)(\CI)$ is reflexive.   
\end{proof}

\subsection{Objects with a Verma flag}\label{Verma flags}
Let $\CG$ be a GKM-graph and choose $x\in\CG$. Let $V(x)\in\CZ\catmod^f$ be the free $S$-module of rank one, generated in degree zero, on which $\CZ$ acts via the evaluation $\CZ\to S$, $(z_v)\mapsto z_x$. We call  $V(x)$ the {\em Verma module} at  $x$. It is obviously flabby. Considered as a sheaf $\SV(x)$ on the moment graph it is just the free skyscraper of rank one at the vertex $x$. Lemma \ref{ext of Vermas} implies the following.
\begin{lemma} Let $M\in\CZ\catmod^f$. Then $M\in\CV=\CZ\catmod^{Verma}$ if and only if $M^{[x]}$ is graded free over $S$, i.e.~ $M^{[x]}\cong \bigoplus V(x)[k_i]$ for some $k_1,\dots,k_l\in\DZ$, for any $x\in\CG$.
\end{lemma}

\subsection{Tilting duality}\label{tilting duality}
For a moment graph $\CG$ let $t(\CG)$ be the {\em tilted moment graph}, i.e.~ the moment graph with the same sets of vertices and edges, the same labelling, but with the reversed order on the vertices. Since the structure algebra does not depend on the order we have
$\CZ(\CG)=\CZ(t(\CG))$.

We can view any $M\in\CV(\CG)$ as a $\CZ$-module. The graded dual $M^\ast=\Hom^\bullet_S(M,S)$ inherits a $\CZ$-structure. This duality reverses the order filtration and since every quotient is free as an $S$-module we get the {\rm tilting equivalence}
$$
t\colon \CV(\CG)\xrightarrow{\sim}\CV(t(\CG))^{opp}
$$
of exact categories (for exactness use Lemma \ref{ext of Vermas}). It corresponds to the tilting equivalence on category $\CO$ for moment graphs associated to non-critical blocks (cf.~ Section \ref{representation theory}).

\section{Projective Objects}
Let $\CA$ be an exact category. Recall that an object $P$ of $\CA$ is called {\em projective}, if the functor
$\Hom_{\CA}(P,\cdot)$ transforms short exact sequences in $\CA$ to short exact sequences of abelian groups.

\subsection{Local characterization}
Let $\CG$ be a GKM-graph. In this section we consider the category $\CC^{ref}$ of reflexive objects and we view any $\SM\in\CC^{ref}$ as a sheaf on $\CG$. Let $\CI\subset\CG$ be open. Recall that we defined the sheaf $\SM^\CI$ by setting $\SM^x$ and $\SM^E$ to zero if $x$ or $E$ do not belong to $\CI$. Since $\SM$ is flabby, $\SM^\CI=\Loc(\Gamma(\SM)^\CI)$, hence our notations are compatible.

\begin{proposition}\label{flabby projectives} Let $\SP\in \CC^{ref}$. Suppose that
\begin{enumerate}
\item for any vertex $x$ of $\CG$ the $S$-module $\SP^x$ is projective and
\item for any edge $E\colon y\xrightarrow{\alpha} x$ the map $\rho_{y,E}\colon \SP^y\to \SP^E$ induces an isomorphism
$\SP^y/\alpha \SP^y \xrightarrow{\sim}  \SP^E$.
\end{enumerate}
Then $\SP$ is projective in $\CC^{ref}$.  
\end{proposition}

\begin{proof}
We have to show that for an admissible epimorphism $\SB \to \SC$ in $\CC^{ref}$ the induced map $\Hom(\SP, \SB)\to\Hom(\SP,\SC)$ is surjective. We first prove this in the case of the quotient $\SM^{\varle x}\to\SM^{<x}$ for $\SM\in\CC^{ref}$ and $x\in\CG$.

Choose a map $f\colon \SP\to \SM^{<x}$. We want to extend it to the vertex $x$, i.e.~ we have to define compatible maps $f^E\colon \SP^E\to \SM^E$ for $E\in\CE_{\delta x}$ and $f^x\colon \SP^x\to \SM^x$. Let $E\colon y\xrightarrow{\alpha} x$ be an edge in $\CE_{\delta x}$. Because of property (2) of $\SP$ and since $\SM^E$ is annihilated by $\alpha$, the map $f^y\colon \SP^y\to \SM^y\xrightarrow{\rho_{y,E}} \SM^E$ determines $f^E\colon \SP^E\to \SM^E$. 

Moreover, $f$ induces a map $f^{<x}\colon\SP(\{<x\})\to \SM(\{<x\})$ and

\centerline{
\xymatrix{
\SP(\{<x\}) \ar[r]\ar[d] & \SM(\{<x\})\ar[d] \\
\bigoplus_{E\in\CE_{\delta x}} \SP^E \ar[r]^{\bigoplus f^E} & \bigoplus_{E\in\CE_{\delta x}} \SM^E
}
}\noindent
commutes. We get an induced map $f^{\delta x}\colon \SP^{\delta x}\to \SM^{\delta x}$. Since $\SM$ is reflexive, it is flabby, so by Proposition \ref{flabby sheaves} the map $\SM^x\to \SM^{\delta x}$ is surjective. Since $\SP^x$ is projective, we can find a map $f^x\colon \SP^x\to \SM^x$ such that

\centerline{
\xymatrix{
\SP^x \ar[d]\ar[r]^{f^x} & \SM^x \ar[d]\\
\SP^{\delta x} \ar[r]^{f^{\delta x}} & \SM^{\delta x} 
}
}
\noindent
commutes. That means that our choices for $f^E$ and $f^x$ are compatible with $\rho$ and we arrive at an extension $\SP^{\leq x}\to \SM^{\leq x}$ of $f$. 

Now let $h\colon\SB\to \SC$ be an admissible epimorphism in $\CC^{ref}$ and choose a map $f\colon \SP\to \SC$. We want to construct $g\colon \SP\to \SB$ such that $f=h\circ g$. We  construct $g$ vertex by vertex, following the order. Suppose we have constructed $g^y$ for all vertices $y$ with  $y<x$. By what we showed above we can find an extension $j^{\leq x}\colon \SP^{\leq x}\to \SB^{\leq x}$ of the map $g^{<x}\colon \SP^{<x}\to \SB^{<x}$. Set 
$$
f^\prime:=h^{\leq x}\circ j^{\leq x}- f^{\leq x}\colon \SP^{\leq x}\to \SC^{\leq x}.
$$
Suppose we find a map $g^\prime\colon \SP^{\leq x}\to \SB^{\leq x}$ such that $f^\prime= h^{\leq x}\circ g^\prime$. Then $f^{\leq x}=h^{\varle x}\circ(j^{\leq x}-g^{\prime})$. So it is enough to find $g^\prime$.

Now $f^\prime\colon \SP^{\leq x}\to \SC^{\leq x}$ is such  that $(f^\prime)^{< x}=0$, i.e.~ $f^\prime$ is given by its component $(f^\prime)^x\colon \SP^x \to \SC^x$ and the image of $(f^\prime)^x$ lies inside $\SC^{[x]}=\ker(\SC^x\to \SC^{\delta x})$. Therefore it is enough to find $(g^\prime)^x\colon \SP^x \to \SB^{[x]}$ such that $(f^\prime)^x=h^{[x]}\circ (g^\prime)^x$. By the exactness of the sequence we started with, the map $h^{[x]}\colon \SB^{[x]} \to \SC^{[x]}$ is surjective, and since $\SP^x$ is projective we find $(g^{\prime})^x\colon \SP^x\to \SB^{[x]}$.
\end{proof}
\subsection{Braden-MacPherson sheaves}\label{Braden-MacPherson sheaves}
In the following we assume that the order on the set of vertices is bounded from above, i.e.~ that for any vertex $v$ of $\CG$ the set $\{x\mid x\geq v\}$ is finite.  
\begin{theorem}\label{projective sheaves} For any $v\in \CG$ there exists an object $\SB(v)\in\CC^{ref}$, unique up to isomorphism, with the following properties:
\begin{enumerate}
\item $\SB(v)$ is indecomposable and projective.
\item $\supp(\SB(v))\subset\{x\mid x\geq v\}$ and $\SB(v)^v\cong S$.
\end{enumerate}
Moreover, every indecomposable projective object in $\CC^{ref}$ is isomorphic to $\SB(v)$ for a unique $v\in\CG$. 
\end{theorem}

\begin{proof}
Let $v$ be a vertex of $\CG$. We first use the Braden-MacPherson algorithm (cf.~ \cite{BMP}, 1.4) to construct a $\CG$-sheaf $\SB(v)$ and then we show that it has all the properties listed in the theorem. Its uniqueness follows from the construction (cf.~ Lemma \ref{automorphisms}).

We  build $\SB(v)$ inductively, using the order on the vertices. 
Suppose we have already constructed $\SB(v)^y$ on a vertex $y$. As Proposition \ref{flabby projectives} suggests, we then define $\SB(v)^E:=\SB(v)^y/\alpha \SB(v)^y$ for any edge $E\colon y\xrightarrow{\alpha} x$ and we let $\rho_{y,E}\colon \SB(v)^y\to \SB(v)^E$ be the canonical quotient map. So we only have to define $\SB(v)$ on the vertices together with the maps $\rho$ in the downwards direction. 

We set $\SB(v)^x=0$ for all edges $x$ with $x\not\geq v$. Then we define $\SB(v)^v=S$. Suppose we have already constructed $\SB(v)^y$ for any vertex $y$ with $y< x$.  As before we define $\SB(v)^{\delta x}\subset \bigoplus_{E\in\CE_{\delta x}} \SB(v)^E$  as the image of $\SB(v)^{<x}(\{<x\})$ and we let $\SB(v)^x$ be a projective (i.e.~ free) cover of $\SB(v)^{\delta x}$ as an $S$-module. The components of the projective cover map $\SB(v)^x\to \SB(v)^{\delta x}\subset \bigoplus \SB(v)^E$ define the maps $\rho_{x,E}\colon \SB(v)^x\to \SB(v)^E$. 

So we have constructed a $\CG$-sheaf $\SB(v)$. By construction $\supp(\SB(v))\subset\{x\mid x\geq v\}$ and $\SB(v)^v=S$. We now show that $\SB(v)$ is generated by global sections, i.e.~ that $(\Loc\circ\Gamma)(\SB(v))\to \SB(v)$ is an isomorphism. Every $m\in \SB(v)^x$ is a component of a global section, hence $\Gamma(\SB(v))^x=\SB(v)^x$ for any vertex $x$. For an edge $E\colon x\xrightarrow{\alpha} y$ consider the commutative diagram

\centerline{
\xymatrix{
(\Loc\circ \Gamma)(\SB(v))^x=\SB(v)^x\ar[d]\ar[dr] & \\
(\Loc\circ\Gamma)(\SB(v))^E \ar[r]^{g^E} & \SB(v)^E.
}
}
\noindent
The vertical and diagonal maps are surjective and $\alpha\SB(v)^x$ is the kernel of the diagonal map. It is contained in the kernel of the vertical map, hence the horizontal map is injective, hence an isomorphism. So $\SB(v)$ is generated by global sections. It is easy to see that $\SB(v)$ is indecomposable. Any $\SB(v)^x$ is free, in particular reflexive, hence Proposition \ref{sections on reflexives} gives $\SB(v)\in\CC^{ref}$.  By Proposition \ref{flabby projectives}, $\SB(v)$ is projective.

The following Lemma follows from the fact that, if $B$ is an $S$-module and $A\to B$ is a projective cover, then any endomorphism on $A$ that induces an automorphism on $B$ is an automorphism. In \cite{BMP}, Proposition 1.2, one can find a stronger statement. 
\begin{lemma}\label{automorphisms}
Let $\SB(v)$ be the $\CG$-sheaf constructed above and let $f\in\End(\SB(v))$ be an endomorphism such that $f^v\in\End_S(\SB(v)^v)$ is an automorphism. Then $f$ is an automorphism.  
\end{lemma}

 Let $\SP\in\CC^{ref}$ be projective and choose $v$ minimal such that $\SP^v\neq 0$. Then $\SP^v$ is free over $S$. Choose a direct summand of rank one in $\SP^v$ and identify it with $\SB(v)^v$. By projectivity the corresponding inclusion and projection are induced by maps $f\colon \SB(v)\to \SP$ and $g\colon \SP\to \SB(v)$. Then $(g\circ f)^v\in\End(\SB(v)^v)$ is the identity, so, by the lemma above, $g\circ f$ is an automorphism. Thus $\SB(v)$ is a direct summand of $\SP$. This proves uniqueness of $\SB(v)$ up to isomorphism and shows that every indecomposable projective in $\CC^{ref}$ is isomorphic to $\SB(v)$ for a unique $v\in\CG$. 
\end{proof}

\subsection{A question}\label{proj and Verma flag} {\em When does $\SB(v)$ admit a Verma flag?}
On the level of $\CG$-sheaves, admitting a Verma flag means that the kernel of the map $\SB(v)^x\to \SB(v)^{\delta x}$ is graded free for all $x\in\CG$. From this perspective this seems to be an interesting and strong additional property of the moment graph. It also implies that the global section of $\SB(v)$ form a free $S$-module.

The results in the following show that the moment graphs associated to blocks outside the critical hyperplanes of the category $\CO$ of a symmetrizable Kac-Moody algebra have this property. The results in \cite{BMP} show that this is also true for moment graphs associated to equivariantly formal varieties.

\subsection{Dependance on the base} Let $R$ be the localization of $S=S(V)$ at the point $0\in V^\ast$, i.e.~ at the maximal ideal $S\cdot V$. For everything we have done so far we can replace $S$ by $R$. We denote the resulting objects by the index $R$. For example, $\CV_R(\CG)$ is the $R$-linear category of objects that admit a Verma flag.

\section{Equivariant Representation Theory}\label{representation theory}
Let $\fg$ be a complex symmetrizable Kac-Moody algebra, $\fb\subset\fg$ its Borel subalgebra and $\fh\subset\fb$ its Cartan subalgebra. Let $U(\fg)$, $U(\fb)$ and $S=S(\fh)=U(\fh)$ be the universal enveloping algebras. We give a short overview of the equivariant representation theory of $\fg$. More details and proofs can be found in \cite{Fie1} and \cite{Fie2}. 

\subsection{Equivariant modules with a Verma flag}
 Let $R=S_{(\fh)}$ be the localization of the symmetric algebra of $\fh$ at the point $0\in\fhd$. We will work with $R$ instead of $S$ since, for most of the following representation-theoretic constructions, it is essential that the base ring is local.

Let $\tau\colon S\to R$ be the canonical map. For $\lambda\in\fhd$ we let $R_\lambda$ be the $\fh$-R-bimodule that is free of rank one over $R$ and on which $\fh$ acts according to the character $\lambda+\tau$ (that means that $H\in\fh$ acts by multiplication with $\lambda(H)\cdot 1+\tau(H)\in R$). We have $\fb=[\fb,\fb]\oplus \fh$ and the projection $\fb\to \fh$ onto the second summand is a morphism of Lie algebras. So
we can consider $R_\lambda$ as a $\fb$-$R$-bimodule and induction yields a $\fg$-$R$-bimodule
$$
M_R(\lambda) = U(\fg)\otimes_{U(\fb)} R_\lambda.
$$
It is called the {\em equivariant} (or {\em deformed}) {\em Verma module with highest weight $\lambda$}. 

Let $M$ be a $\fg$-$R$-bimodule and $\mu\in\fhd$. Then
$$
M_\mu:=\left\{m\in M\mid H.m=(\mu(H)\cdot 1+\tau(H))m\quad\text{for all $H\in\fh$}\right\}
$$
is called its {\em weight space of weight $\mu$}. We say that $M$ is an  {\em (equivariant) weight module} if $M=\bigoplus_{\mu\in\fhd} M_\mu$ and we call $\{\mu\in\fhd\mid M_\mu\neq0\}$ its set of weights. The equivariant Verma module $M_R(\lambda)$ is a weight module. Its set of weights is $\lambda-\DN\cdot \Delta^{+}$, where $\Delta^+\subset\fhd$ is the set of roots of $\fb$ (i.e.~ the set of positive roots).

Let $\CM_R$ be the category of $\fg$-$R$-bimodules that are weight modules and that admit a finite Verma flag, i.e.~ a finite filtration whose subquotients are isomorphic to equivariant Verma modules. Let  $\left(M:M_R(\lambda)\right)$ be the multiplicity of $M_R(\lambda)$. It is independent of the filtration.
$\CM_R$ is a subcategory of the category $\CO_R$, which is the category of $\fg$-$R$-bimodules $M$ that are weight modules and locally finite for the action of $\fb$, i.e.~ that have the property that for any $m\in M$ the space $(U(\fb)\otimes_\DC R).m$ is finitely generated over $R$. Note that $\CO_R$ is abelian and $\CM_R$ is closed under extensions in $\CO_R$, so $\CM_R$ inherits an exact structure.

\subsection{Projective objects}
Let $\Delta\subset\fhd$ be the set of roots of $\fg$ and $\Delta_+\subset\Delta$ the set of roots of $\fb$. We define a partial order on $\fhd$ by setting $\lambda\geq\mu$ if and only if $\lambda-\mu$ can be written as a sum of positive roots.  

For $\nu\in\fhd$ let $\CM_R^{\varle \nu}$ be the full subcategory of $\CM_R$ consisting of modules $M$ such that $M_\mu\neq0$ implies $\mu\leq\nu$. This is equivalent to saying that only Verma modules with highest weights $\mu$ with $\mu\leq\nu$  occur in a Verma flag.  $\CM^{\varle\nu}_R$ is a subcategory of the analogously defined category $\CO_R^{\varle\nu}\subset\CO_R$ and it is closed under extensions, so it inherits an exact structure.

The category $\CO_R^{\varle\nu}$ contains enough projective objects and all of them admit a Verma flag, i.e.~ they are projective in $\CM_R^{\varle\nu}$. For any $\lambda$ with $\lambda\leq \nu$ there exists an object $P_R^{\varle \nu}(\lambda)\in\CM_R^{\varle\nu}$ which is characterized by being indecomposable and projective and the properties
\begin{eqnarray*}
\left(P_R^{\varle\nu}(\lambda)\colon M_R(\mu)\right) & \neq & 0 \Rightarrow \mu\geq \lambda \quad\text{(and $\mu\leq\nu$),} \\
\left(P_R^{\varle\nu}(\lambda)\colon M_R(\lambda)\right) & = & 1. 
\end{eqnarray*}
Every indecomposable projective object in $\CO_R^{\varle\nu}$ is isomorphic to  $P_R^{\varle\nu}(\lambda)$ for a unique $\lambda$.

\subsection{Block decomposition}\label{block decomposition}

For a symmetrizable Kac-Moody algebra $\fg$ there exists an invariant, non-degenerate  bilinear form $(\cdot,\cdot)\colon\fg × \fg\to\DC$. It induces a non-degenerate bilinear form $(\cdot,\cdot)\colon \fh × \fh\to \DC$ and, dually, a non-degenerate bilinear form $(\cdot,\cdot)\colon \fhd × \fhd\to \DC$. 
Choose $\rho\in\fhd$ such that $(\rho,\alpha)=1$ for any simple root $\alpha\in\Delta_+$.

Let $\sim$ be the finest equivalence relation on $\fhd$ such that $\lambda\sim\mu$ if there exist $\beta\in\Delta$ and $n\in\DZ$ with $\lambda-\mu=n\beta$ and $2(\lambda+\rho,\beta)=n(\beta,\beta)$. This definition does not depend on the choice of $\rho$.
For an equivalence class $\Lambda\in\fhd/_{\textstyle{\sim}}$ let $\CM_{R,\Lambda}$ be the subcategory of $\CM_R$ consisting of all modules $M$ such that $\left(M:M_R(\lambda)\right)\neq 0$ implies $\lambda\in\Lambda$. The next result follows from the analogous result in \cite{Fie1}, Proposition 2.8.

\begin{theorem}
The functor
\begin{eqnarray*}
\prod_{\Lambda\in\fhd/_{\scriptstyle{\sim}}} \CM_{R,\Lambda} & \to  & \CM_R, \\
\{M_\Lambda\}   & \mapsto  &    \bigoplus_{\Lambda\in\fhd/_{\scriptstyle{\sim}}} M_\Lambda
\end{eqnarray*}
 is an equivalence of categories.
\end{theorem}

Let $\Lambda\in\fhd/_{{\textstyle{\sim}}}$ be an equivalence class. For $\nu\in\Lambda$ let $\Lambda^{\varle\nu}\subset\Lambda$ be the subset of elements that are smaller or equal to $\nu$. Let $\CM_{R,\Lambda^{\varle\nu}}=\CM_R^{\varle\nu}\cap \CM_{R,\Lambda}$ and $\CO_{R,\Lambda^{\varle\nu}}=\CO_R^{\varle\nu}\cap\CO_{R,\Lambda}$ be the respective truncations. Then $M\in\CM_R$ lies in $\CM_{R,\Lambda^{\varle\nu}}$ if and only if $\left(M:M_R(\mu)\right)\neq 0$ implies $\mu\in\Lambda^{\varle\nu}$. Let $\CP_{R,\Lambda^{\varle\nu}}\subset\CM_{R,\Lambda^{\varle\nu}}$ be the full subcategory of projective objects. The set $\{P_R^{\varle\nu}(\lambda)\}_{\lambda\in\Lambda^{\varle\nu}}$ is a faithful set of small projectives in the sense of \cite{Mit}. If we denote by $R\catmod^{(\CP_{R,\Lambda^{\varle\nu}})^{opp}}$ the category of all additive contravariant functors from $\CP_{R,\Lambda^{\varle\nu}}$ to $R\catmod$, then the functor
\begin{eqnarray*}
\CO_{R,\Lambda^{\varle\nu}} & \to & R\catmod^{(\CP_{R,\Lambda^{\varle\nu}})^{opp}} \\
M & \mapsto & M(\cdot):=\Hom_{\CO_R}(\cdot, M)
\end{eqnarray*}
is an equivalence of categories. 

\subsection{Canonical Verma flags}
Let $\Lambda\in\fhd/_{\textstyle{\sim}}$ and $M\in\CO_{R,\Lambda}$. For any $\mu\in\Lambda$ let $M_{\varge \mu}\subset M$ be the sub-$\fg$-$R$-bimodule  that is generated by all weight spaces $M_\nu$ with $\nu\geq\mu$ and let $M^{\varle \mu}$ be the cokernel of $\sum_{\nu\not\varle \mu} M_{\varge \nu}\subset M$. Define $M^{<\mu}$ analogously. For $\mu,\mu^\prime\in\Lambda$ with $\mu\geq \mu^\prime$ we have a canonical surjection $M^{\varle \mu}\to M^{\varle \mu^\prime}$, i.e.~ $M^{\varle \cdot}$ is a cofiltration indexed by the ordered set $\Lambda^{opp}$.  $M$ admits a Verma flag if and only if for any $\mu\in\Lambda$ the kernel $M^{[\mu]}:=\ker(M^{\varle \mu}\to M^{<\mu})$ is isomorphic to a finite direct sum of copies of the Verma module  $M_R(\mu)$ and only finitely many of them are non-zero. If $M$ admits a Verma flag 
we call $M^{\varle\cdot}$ the {\em canonical Verma cofiltration} of $M$.

\subsection{The associated integral Coxeter system} Let $\CW$ be the Weyl group of $\fg$. It naturally acts on $\fhd$. Recall that a root $\alpha\in\Delta$ is called {\em real} if it is $\CW$-conjugate to a simple root. Roots that are not real are called {\em imaginary}. Let $\Delta^{re}$ be the set of real roots and $\Delta^{im}$ the set of imaginary roots. For any $\beta\in\Delta^{re}$ let $s_\beta\in\CW$ be the corresponding reflection.

There is an elegant description of the non-critical equivalence classes $\Lambda\in\fhd/_{\textstyle{\sim}}$ as orbits of certain subgroups of $\CW$. We say that an element $\lambda\in\fhd$ lies {\em outside the critical hyperplanes}, if $2(\lambda+\rho,\beta)\not\in \DZ(\beta,\beta)$ for any imaginary root $\beta\in\Delta^{im}$. We say that an equivalence class $\Lambda\in\fhd/_{\textstyle{\sim}}$ lies {\em outside the critical hyperplanes}, if any $\lambda\in\Lambda$ does. For such $\Lambda$ let 
\begin{eqnarray*}
\Delta_\Lambda & := & \left\{\beta\in\Delta\mid 2(\lambda+\rho,\beta)\in\DZ(\beta,\beta)\text{\ for some  $\lambda\in\Lambda$}\right\}, \\
\CW_\Lambda & := & \langle s_\beta\mid \beta\in\Delta_\Lambda\rangle\subset\CW, \\  
\Delta_{\Lambda,+}& :=& \Delta_\Lambda\cap\Delta_+, \\ 
\Pi_\Lambda & := & \{\alpha\in\Delta_{\Lambda,+}\mid s_\alpha(\Delta_{\Lambda,+}\setminus\{\alpha\})\subset\Delta_{\Lambda,+}\}, \\ 
\CS_\Lambda&  := & \{s_\alpha\mid\alpha\in\Pi_\Lambda\}.
\end{eqnarray*}
Then $(\CW_\Lambda,\CS_\Lambda)$ is a Coxeter system (it is called the {\em integral Coxeter system for $\Lambda$}). Denote the linear $\CW$-action on $\fhd$ by $\lambda\mapsto w(\lambda)$ and let $\lambda\mapsto w.\lambda=w(\lambda+\rho)-\rho$ be the action shifted by $\rho$. Again it does not depend on the choice of $\rho$. Then $\Lambda$ is a $\CW_\Lambda$-orbit in $\fhd$, i.e.~ $\Lambda=\CW_\Lambda.\lambda$
for any $\lambda\in\Lambda$. 

\subsection{The associated moment graph}\label{associated moment graph}
Let $\Lambda\in\fhd/_{\textstyle{\sim}}$ be an equivalence class outside the critical hyperplanes. The {\em associated moment graph} $\CG=\CG_\Lambda$ over $\fh$ is defined as follows.  Its  set of vertices is $\Lambda$. Let $\CT_\Lambda\subset\CW_\Lambda$ be the set of reflections in $\CW_{\Lambda}$ and for $t\in\CT_{\Lambda}$ let $\alpha^\lor_t\in\fh$ be the corresponding coroot.  The vertices $\lambda,\mu\in\Lambda$ are linked by an edge if there is a reflection $t\in\CT_{\Lambda}$ such that $\lambda=t.\mu$ and $\lambda\neq\mu$. We label the edge by $\DC\cdot\alpha^\lor_t\subset\fh$. The order on the vertices is the induced order on $\Lambda\subset \fhd$. 

If $\lambda\in\Lambda$ is antidominant, i.e.~ a smallest element in $\Lambda$, then $\CW_\Lambda/\CW_{\Lambda}^\lambda\to\Lambda, \bar w\mapsto \bar w.\lambda$ is a bijection of ordered sets and $\CG_\Lambda$ is isomorphic to the moment graph $\CG(\CW_{\Lambda},\CS_{\Lambda},\CW_{\Lambda}^\lambda)$ defined in  \ref{examples of moment graphs}. If $\lambda\in\Lambda$ is dominant, i.e.~ a largest element in $\Lambda$ with respect to the Bruhat order, then $\CG_\Lambda$ is isomorphic to the tilted moment graph 
$t(\CG(\CW_{\Lambda},\CS_{\Lambda},\CW_\Lambda^\lambda))$. By Lemma \ref{quasi-finiteness of graphs}, $\CG_\Lambda$ is quasi-finite and it has the GKM-property since any two coroots that induce different reflections are linearly independent.  

\subsection{The equivariant center}
Let $\Lambda\in\fhd/_{\textstyle{\sim}}$ be a non-critical equivalence class. Let $\CZ_{R,\Lambda}$  be the center of $\CO_{R,\Lambda}$, i.e.~ the ring $\End(\id_{\CO_{R,\Lambda}})$ of endotransformations of the identity functor on $\CO_{R,\Lambda}$.  It is naturally isomorphic to the center of $\CM_{R,\Lambda}$ since all projectives in $\CO_{R,\Lambda}$ lie in $\CM_{R,\Lambda}$. Let 
$$
ev\colon \CZ_{R,\Lambda}\to \prod_{\nu\in\Lambda} \End(M_R(\nu))= \prod_{\nu\in\Lambda} R
$$
be the evaluation at the Verma modules. 
\begin{theorem}[\cite{Fie1}, Theorem 3.6]\label{calculation of the center} The evaluation map is injective and defines an isomorphism
$$
\CZ_{R,\Lambda} \cong \left\{(z_{\nu})\in\prod_{\nu\in\Lambda} R
\left|
\,
\begin{matrix}
 z_{\nu}\equiv z_{t.\nu} \mod \alpha^\lor_t \\
\text{for all $\nu\in\Lambda$, $t\in\CT_{\Lambda}$} 
\end{matrix}
\right.
\right\}.
$$
\end{theorem}
Hence the center of a block $\CO_{R,\Lambda}$ is canonically isomorphic to the structure algebra $\CZ_R(\CG_\Lambda)$ of its associated moment graph. Let $\nu\in\Lambda$ and let $\CZ_{R,\Lambda^{\varle \nu}}$ be the center of the truncated category $\CO_{R,\Lambda^{\varle\nu}}$. Analogously we have an isomorphism
$$
\CZ_{R,\Lambda^{\varle \nu}} = \left\{
\left.
(z_{\mu})\in\prod_{\mu\in\Lambda^{\varle\nu}} R
\,
\right|
\begin{matrix}
 z_{\mu}\equiv z_{t.\mu} \mod \alpha^\lor_t \\
\text{for all $\mu,t.\mu\in\Lambda^{\varle\nu}$, $t\in\CT_{\Lambda}$}
\end{matrix}
\right\}.
$$
So if $\CG_\Lambda^{\varle\nu}$ is the restricted moment graph, i.e.~ the subgraph that contains only the vertices that are smaller or equal to $\nu$, we again have an identification $\CZ_{R,\Lambda^{\varle\nu}}=\CZ_R(\CG_\Lambda^{\varle\nu})$. 

\begin{remark}
It is a non-trivial fact that the restriction $\CZ_{R,\Lambda}\to\CZ_{R,\Lambda^{\varle\nu}}$ is surjective. This corresponds to the smoothness of flag varieties and to the fact that the multiplicities in the antidominant projective are at most one (cf.~ \cite{Fie2}, Lemma 3.8). In the language of moment graphs this means that the structure sheaf $\SA_{R,\CG_\Lambda}$ is flabby.
\end{remark}
In the following we denote by $\CV_{R,\Lambda}$ and $\CV_{R,\Lambda^{\varle\nu}}$ the categories $\CV_R(\CG_\Lambda)$ and $\CV_R(\CG_\Lambda^{\varle\nu})$.  

\section{An equivalence of categories}
In this section we state and prove the main result, namely Theorem \ref{equivalence of categories}.

\subsection{The structure functor}\label{structure functor}

Let $\Lambda\subset\fhd/_{\textstyle{\sim}}$ be an equivalence class outside the critical hyperplanes. In \cite{Fie2} we defined the  {\em structure functor}
$$
\DV=\DV_{R,\Lambda}\colon \CM_{R,\Lambda} \to \CZ_{R,\Lambda}\catmod^f.
$$
If $\Lambda$ contains an antidominant, i.e.~ smallest  element $\lambda$, then $\DV=\Hom(P_R(\lambda),\cdot)$, where $P_R(\lambda)=\varprojlim P^{\varle \nu}_R(\lambda)$ is a completion of the truncated antidominant projective objects. If $\Lambda$ contains a dominant, i.e.~ largest element we use the tilting functor to translate the situation to the antidominant case. 

Note that $\DV$ is exact for the classical exact structure on $\CZ_{R,\Lambda}\catmod^f$ since $P_R(\lambda)$ is projective.
Let $M\in\CM_{R,\Lambda}$. Then $M\otimes_R Q\in\CM_{Q,\Lambda}$ splits into the direct sum of Verma modules. Applying the structure functor gives the canoncial decomposition of $\DV M\otimes_R Q$ in \ref{decomposition at generic point}. Hence $\DV$ transforms the canonical Verma filtration  to the order filtration. In particular, $\DV$  is exact even for the non-standard exact structure defined in \ref{subsection - exact structure}. 

By \cite{Fie2}, {Proposition 3.11}, the image of a Verma module is a Verma module, i.e.~ $\DV M_R(\lambda)\cong \SV_R(\lambda)$, hence the image of $\CM_{R,\Lambda}$ is contained in $\CV_{R,\Lambda}$. We denote the restriction $\DV\colon \CM_{R,\Lambda}\to \CV_{R,\Lambda}$ again by $\DV$. 

\begin{theorem} \label{equivalence of categories} The functor 
$\DV\colon  \CM_{R,\Lambda}\to \CV_{R,\Lambda}$ 
is an equivalence of exact categories. 
\end{theorem}

\subsection{The structure of projectives}
For $\nu\in\Lambda$ the structure functor restricts to a functor $\DV\colon\CM_{R,\Lambda^{\varle\nu}}\to\CV_{R,\Lambda^{\varle\nu}}$ that we denote by $\DV$ as well. 
The most important step in the proof of Theorem \ref{equivalence of categories} is the following proposition.
\begin{proposition}\label{v of projectives} Let $P\in\CM_{R,\Lambda^{\varle \nu}}$ be projective. Then $\DV P$ is projective in $\CV_{R,\Lambda^{\varle \nu}}$.
\end{proposition}

\begin{proof} We will show that $\DV P$ has the two properties listed in Proposition \ref{flabby projectives}. In the following we make extensive use of the base change results that are described in \cite{Fie1}, 2.3, and \cite{Fie2}.

Let $\chi\in\Lambda^{\varle\nu}$. We want to show that $(\DV P)^{\chi}$ is free over $R$. Let $P_\DC:=P\otimes_R \DC$ and $M_\DC(\chi)=M_R(\chi)\otimes_R \DC$. Then $P_\DC$ is  projective in the classical category $\CO$ and $M_{\DC}(\chi)$ is a classical Verma module. Moreover
$$
\left(P:M_R(\chi)\right)=\left(P_\DC:M_\DC(\chi)\right)=\dim_\DC\Hom_{\CO}(P_\DC,M_{\DC}(\chi)).
$$
Denote this number by $k$.
 
We have $\Hom_{\CO_R}(P,M_R(\chi))\otimes_R \DC=\Hom_{\CO}(P_\DC,M_{\DC}(\chi))$. So we can choose  $f_1,\dots, f_k\in\Hom_{\CO_R}(P,M_R(\chi))$ such that $\bar f_1,\dots,\bar f_k$ with $\bar f_i=f_i\otimes 1$ is a basis of $\Hom_{\CO}(P_\DC,M_{\DC}(\chi))$.  Consider 
$$
f=(f_1,\dots, f_k)\colon P\to  M_R(\chi)^{\bigoplus k}=\bigoplus_{1\leq i\leq k}M_R(\chi).
$$ 
The map $\DV f\colon\DV P\to\DV M_R(\chi)^{\oplus k}$ factors over $\DV P \to (\DV P)^\chi$, since  $\DV M_R(\chi)$ is supported on $\chi$. We claim that 
$(\DV f)^\chi\colon(\DV P)^\chi\to \DV M_R(\chi)^{\bigoplus k}$ 
is an isomorphism of $R$-modules. As $\DV M_R(\chi)$ is free over $R$ this shows the first property in Proposition \ref{flabby projectives}. 

Let $Q$ be the quotient field of $R$. 
For an $R$-module $M$ we write $M_Q$ for its localization $M\otimes_R Q$. We have a canonical decomposition   
$P_Q=\bigoplus_{\mu\in\Lambda^{\varle\nu}} P_{Q,\mu}$
such that each $P_{Q,\mu}$ is isomorphic to a direct sum of copies of $M_Q(\mu)$. Then $f$ induces an isomorphism $P_{Q,\chi}\cong M_Q(\chi)^{\bigoplus k}$, hence $\DV(f_Q)^{\chi}$ is an isomorphism. As $\DV$ and $(\cdot)^\chi$ commute with base change, $(\DV f)^\chi_Q$ is an isomorphism. Since $(\DV P)^{\chi}$ is torsion free, $(\DV f)^{\chi}$ is injective. So we are left with showing that $(\DV f)^{\chi}$ is surjective. 

It is enough to show that $\DV f$ is surjective. It suffices to show that it is surjective at the closed point $\Spec(\DC)$, i.e.~ we have to show that for $\bar f=(\bar f_1,\dots, \bar f_k)\colon P_\DC\to M_\DC(\chi)^{\bigoplus k}$ the map 
$\DV \bar f\colon \DV P_\DC\to \DV M_\DC(\chi)^{\bigoplus k}$ 
is surjective.

We have to distinguish two cases. First we assume that $\Lambda$ contains an antidominant element $\lambda$. 
 Then $\DV \cdot =\Hom_{\CO}(P_\DC(\lambda), \cdot)$. Let
$P_\DC^{\chi}\subset M_\DC(\chi)^{\bigoplus k}$ be the image of $\bar f$. By projectivity of $P_\DC(\lambda)$, $\im\,\DV \bar f=\DV P_\DC^\chi$. Since $\DV M_\DC(\chi)$ is one-dimensional, we have to show that the dimension of $\DV P_\DC^\chi$ is at least $k$. 

\begin{lemma} The socle of $M_\DC(\chi)$ is simple and isomorphic to $M_\DC(\lambda)$. 
\end{lemma}
Before we prove the lemma, we finish our previous argument.
$\DV M_\DC(\chi)$ is generated by the composition $P_\DC(\lambda)\to M_\DC(\lambda)\inj M_\DC(\chi)$. Let $g\colon P_\DC^{\chi}\to M_\DC(\chi)$ be a non-zero map. By the lemma, $M_\DC(\lambda)$ is contained in the image of $g$ and by projectivity of $P_\DC(\lambda)$ any map $P_\DC(\lambda)\to M_\DC(\chi)$ factors over $g$. In particular, the map $\Hom(P_{\DC}^\chi, M_\DC(\chi))\to \Hom(\DV P_{\DC}^\chi,\DV M_\DC(\chi))$, that is induced by $\DV$, is injective.

Every homomorphism $P_\DC\to M_\DC(\chi)$ factors over $P_\DC^\chi$, and the dimension of $\Hom(P_\DC,M_\DC(\chi))$ is $k$. Hence $\dim_\DC \Hom(\DV P_{\DC}^\chi,\DV M_\DC(\chi))\geq k$. Since $\DV M_\DC(\chi)$ is one-dimensional, the dimension of $\DV P_{\DC}^\chi$ is at least $k$. This is what we wanted to show. 

\begin{proof}[Proof of the Lemma] 
We use induction on the Bruhat order on $\Lambda$. The case $\chi=\lambda$ follows from the fact that $M_\DC(\lambda)$ is simple. Let $\chi\in\Lambda$ and let $M\subset M_\DC(\chi)$ be a non-trivial submodule. Let $\nu<\chi$ be a maximal weight of $M$. Then there is a non-zero map $i\colon M_\DC(\nu)\to M$. Every non-trivial map between Verma modules is injective, so by the induction hypothesis there is a non-zero map $M_\DC(\lambda)\to M$. So any submodule of $M_\DC(\chi)$ contains a copy  of $M_\DC(\lambda)$. Since the multiplicity of the simple module $M_\DC(\lambda)$ in $M_\DC(\chi)$ is one (cf.~ \cite{Fie2}, Lemma 3.8),  this is the unique simple submodule. 
\end{proof}

Now assume that $\Lambda$ contains a largest element $\lambda$. Recall that we want to show that the map $\DV \bar f\colon \DV P_\DC \to \DV M_\DC(\chi)^{\bigoplus k}$ is surjective. Let $\CM\subset\CO$ be the non-equivariant category of modules that admit a Verma flag and $\CM_\Omega=\CM\cap\CO_\Omega$ the corresponding blocks.
There is a {\em tilting equivalence} $t\colon \CM_{\Lambda}\to \CM_{t(\Lambda)}^{opp}$, where $t(\Lambda):=-2\rho-\Lambda$ is an equivalence class with smallest element $-2\rho-\lambda$ (cf.~ \cite{Soe2} and \cite{Fie2}). For $M\in\CM_\Lambda$ we have $\DV M=\Hom(P_{\DC}(-2\rho-\lambda), t(M))^{\star}$, where $\star=\Hom(\cdot, \DC)$ is the duality functor on $\CZ_{\DC,\Lambda}\catmod$.

The next lemma shows that $t(\bar f)\colon t(M_\DC(\chi)^{\bigoplus k})\to t(P_\DC)$ is injective. Hence $\Hom(P_{\DC}(-2\rho-\lambda), t(M_\DC(\chi)^{\bigoplus k}))\to\Hom(P_{\DC}(-2\rho-\lambda), t(P_\DC))$ is injective, so the dual map $\DV P_\DC\to\DV M_\DC(\chi)^{\bigoplus k}$ is surjective. 

\begin{lemma} Choose $M\in\CM$, $\mu\in\fhd$ and a map $g\colon M_\DC(\mu)^{\oplus k}\to M$. Suppose that $g_\mu\colon  M_\DC(\mu)^{\oplus k}_\mu\to M_\mu$ is injective. Then $g$ is injective.
\end{lemma}
\begin{proof}
Set $N=M_\DC(\mu)^{\oplus k}$. Choose a filtration $M=M_l\supset M_{l-1}\supset\dots\supset M_0=0$ of $M$ by Verma modules. It induces, via $g_\mu$, a filtration $N_\mu=N_{\mu,l}\supset\dots\supset N_{\mu,0}=0$ on the $\mu$-weight space of $N$. Let $N_i\subset N$ be generated by $N_{\mu,i}$. We get a map $g_i\colon N_i/N_{i-1}\to M_i/M_{i-1}$ that is injective on the $\mu$-weight space. Now $M_i/M_{i-1}$ is a Verma module and any non-trivial morphism between Verma modules is injective. Hence $g_i$ is injective. This is the case for any $i$, so $g$ is injective.
\end{proof}

So we proved that $(\DV P)^\chi$ is free over $R$ in both cases.
We still have to show the second property in Proposition \ref{flabby projectives}. Let $E\colon \chi\xrightarrow{{\alpha^\lor}} {\chi^\prime}$ be an edge, starting at $\chi$, i.e.~ $\chi< \chi^\prime$. Since $\DV P$ has a Verma flag, it is determined by local relations, so we can identify $(\DV P)^{\chi}\to (\DV P)^E$ with the quotient map $(\DV P)^{\chi}\to (\DV P)^{\chi}/(\DV P)_{\chi}$. So we have to show that $(\DV P)_{\chi}={\alpha^\lor} (\DV P)^{\chi}$. It is clear that ${\alpha^\lor} (\DV P)^{\chi}\subset (\DV P)_{\chi}$, hence we must show that any $m\in (\DV P)_{\chi}$ is divisible by ${\alpha^\lor}$ in $(\DV P)^{\chi}$.

Denote by $M_{\alpha^\lor}$ the localization $M\otimes_R R_{\alpha^\lor}$ of an $R$-module $M$ at the ideal generated by ${\alpha^\lor}$.  It is enough to show that $m\in (\DV P)_{\chi}$ is divisible by ${\alpha^\lor}$ in the localization $((\DV P)^{\chi})_{\alpha^\lor}$. Note that $((\DV P)^{\chi})_{\alpha^\lor}=(\DV P_{\alpha^\lor})^{\chi}$ since $\DV$ and $(\cdot)^\chi$ commute with base change. 

The module $P_{\alpha^\lor}$ is projective in $\CM_{\alpha^\lor}$, the corresponding version of $\CM$ over $R_{\alpha^\lor}$ (cf.~ \cite{Fie1}, Proposition 2.4). The category $\CM_{\alpha^\lor}$ splits into generic and subgeneric blocks. Since $\chi<\chi^\prime$, the only contribution to $(\DV P_{\alpha^\lor})^{\chi}$ comes from the direct summands isomorphic to the antidominant projective object $P_{\alpha^\lor}({\chi})$. By \cite{Fie1}, Lemma 3.12,  $\DV P_{{\alpha^\lor}}({\chi})$ is isomorphic to the local structure algebra $\CZ_{{\alpha^\lor}}(E)$. But $\CZ_{{\alpha^\lor}}(E)_{\chi}={\alpha^\lor}\cdot \CZ_{{\alpha^\lor}}(E)^{\chi}$, hence $(\DV P_{\alpha^\lor})_{\chi}={\alpha^\lor}\cdot (\DV P_{\alpha^\lor})^{\chi}$. 
\end{proof}

\subsection{The representation functor}

Recall the tilting duality $t\colon \CM_{R,\Lambda}\xrightarrow{\sim}\CM_{R,t(\Lambda)}^{opp}$ (cf.~\cite{Fie2}, 2.6), where $t(\Lambda)=-2\rho-\Lambda$. By definition we have $\DV_\Lambda=\ast\circ\DV_{t(\Lambda)}\circ t$, where $(\cdot)^\ast=\Hom_R(\cdot,R)$.  We can identify $\CG_{t(\Lambda)}$ with the tilted moment graph $t(\CG_\Lambda)$. Using the tilting duality $\CV_{R,\Lambda}\xrightarrow{\sim}\CV_{R,t(\Lambda)}^{opp}$ (cf.~ \ref{tilting duality})
we see that it is enough to consider the case that $\Lambda$ contains an antidominant element in order to prove Theorem \ref{equivalence of categories}. 

Since every object in $\CM_{R,\Lambda}$ or in $\CV_{R,\Lambda}$ lies in some suitably truncated subcategory, we only have to show that $\DV\colon \CM_{R,\Lambda^{\varle \nu}}\to \CV_{R,\Lambda^{\varle\nu}}$ is an equivalence. For this we construct an inverse functor. 

We have an equivalence $\CO_{R,\Lambda^{\varle\nu}}\xrightarrow{\sim}   R\catmod^{(\CP_{R,\Lambda^{\varle\nu}})^{opp}}$
(cf.~ \ref{block decomposition}).
For  $\SM\in\CV_{R,\Lambda^{\varle \nu}}$ and $P\in\CP_{R,\Lambda^{\varle\nu}}$ we set 
$(\DW \SM)(P):=\Hom_{\CV_{R}}(\DV P,\SM)$. 
This yields an additive contravariant functor $\DW \SM\colon \CP_{R,\Lambda^{\varle\nu}} \to R\catmod$, hence an object  $\DW \SM\in\CO_{R,\Lambda^{\varle\nu}}$ and we get a functor 
$$
\DW\colon \CV_{R,\Lambda^{\varle \nu}}\to \CO_{R,\Lambda^{\varle \nu}}.
$$ 
By  Proposition \ref{v of projectives},  $\DV P$ is projective in $\CV_{R,\Lambda^{\varle \nu}}$ for any $P\in\CP_{R,\Lambda^{\varle \nu}}$, hence $\DW$ is an exact functor.

\subsection{Natural transformations}

Let $P\in\CP_{R,\Lambda^{\varle \nu}}$ and $M\in\CM_{R,\Lambda^{\varle \nu}}$. The functor $\DV$ induces a map $\Hom_{\CO_R}(P,M)\to \Hom_{\CV_R}(\DV P,\DV M)$, hence a map $M(P)\to(\DW\circ \DV)(M)(P)$ that is both functorial in $M$ and $P$, hence we get a natural transformation 
$\id_{\CM_{R,\Lambda^{\varle \nu}}}\to\DW\circ\DV$. 

\begin{proposition}[\cite{Fie2}, Theorem 3.15]
For $M,N\in\CM_{R,\Lambda^{\varle \nu}}$ the map
$$
\Hom_{\CO_R}(M,N)\to \Hom_{\CV_R}(\DV  M,\DV N)
$$
is an isomorphism.
\end{proposition}

\begin{proof}[Proof of Theorem \ref{equivalence of categories}]
The proposition above means that the natural transformation 
$\id_{\CM_{R,\Lambda^{\varle\nu}}}\xrightarrow{\sim}\DW\circ\DV$
is an isomorphism. Let $\mu\in\Lambda^{\varle\nu}$. 
From $\DV M_R(\mu)\cong \SV_R(\mu)$ (cf.~ \ref{Verma flags}) we get $\DW \SV_R(\mu)\cong M_R(\mu)$. Since $\DW$ is exact, for any $\SM\in\CV_{R,\Lambda^{\varle \nu}}$ the module $\DW \SM$ has a Verma flag. Hence we get a functor
$$
\DW\colon \CV_{R,\Lambda^{\varle \nu}}\to \CM_{R,\Lambda^{\varle\nu}}.
$$

For any projective $P\in \CP_{R,\Lambda^{\varle \nu}}$, the module $\DV P$ is projective in $\CV_{R,\Lambda^{\varle \nu}}$. If $P$ is indecomposable, then so is $\DV P$. Since both $\DV$ and $\DW$ do not change the support, we get
$\DV P^{\varle\nu}_R(\mu)  \cong  \SB^{\varle\nu}_R(\mu)$ and $\DW \SB^{\varle\nu}_R(\mu) \cong  P^{\varle\nu}_R(\mu)$, where $\SB^{\varle\nu}_R(\mu)$ is the Braden--MacPherson sheaf on $\CG_\Lambda^{\varle\nu}$. 

All that remains to show is that $\DV\circ\DW$ is isomorphic to the identity functor on $\CV_{R,\Lambda^{\varle\nu}}$. By Lemma 3.12 in \cite{Fie2} we have a natural identification $\CZ_{R,\Lambda}^{\varle \nu}= \End(P^{\varle \nu}_R(\lambda))= \DV P^{\varle \nu}_R(\lambda)$. Moreover, for any $\SM\in\CV_{R,\Lambda^{\varle \nu}}$, the object $\DW \SM$ is characterized by $\Hom_{\CO_R}(P,\DW \SM)=\Hom_{\CV_R}(\DV P,\SM)$ for any $P\in\CP_{R,\Lambda^{\varle \nu}}$. Hence  
\begin{eqnarray*}
(\DV\circ\DW)(\SM) & = & \Hom_{\CO_R}(P_R^{\varle\nu}(\lambda), \DW \SM) \\
 & = & \Hom_{\CV_R}(\DV P^{\varle\nu}_R(\lambda),  \SM) \\
 & = & \Hom_{\CZ_R}(\CZ_{R,\Lambda}^{\varle \nu}, \SM ) \\
& = & \SM.
\end{eqnarray*}
This finishes the proof of Theorem \ref{equivalence of categories}. 
\end{proof}

\end{document}